\documentclass[12pt]{article}

\usepackage{amsmath,amssymb,amsthm,amscd,a4wide}

\begin{document}

\newcommand{\End}{{\rm{End}\ts}}
\newcommand{\Hom}{{\rm{Hom}}}
\newcommand{\ch}{{\rm{ch}\ts}}
\newcommand{\non}{\nonumber}
\newcommand{\wt}{\widetilde}
\newcommand{\wh}{\widehat}
\newcommand{\ot}{\otimes}
\newcommand{\la}{\lambda}
\newcommand{\La}{\Lambda}
\newcommand{\De}{\Delta}
\newcommand{\al}{\alpha}
\newcommand{\be}{\beta}
\newcommand{\ga}{\gamma}
\newcommand{\Ga}{\Gamma}
\newcommand{\ep}{\epsilon}
\newcommand{\ka}{\kappa}
\newcommand{\vk}{\varkappa}
\newcommand{\si}{\sigma}
\newcommand{\vp}{\varphi}
\newcommand{\de}{\delta}
\newcommand{\ze}{\zeta}
\newcommand{\om}{\omega}
\newcommand{\ee}{\epsilon^{}}
\newcommand{\su}{s^{}}
\newcommand{\hra}{\hookrightarrow}
\newcommand{\ve}{\varepsilon}
\newcommand{\ts}{\,}
\newcommand{\vac}{\mathbf{1}}
\newcommand{\di}{\partial}
\newcommand{\qin}{q^{-1}}
\newcommand{\tss}{\hspace{1pt}}
\newcommand{\Sr}{ {\rm S}}
\newcommand{\U}{ {\rm U}}
\newcommand{\BL}{ {\overline L}}
\newcommand{\BE}{ {\overline E}}
\newcommand{\BP}{ {\overline P}}
\newcommand{\AAb}{\mathbb{A}\tss}
\newcommand{\CC}{\mathbb{C}\tss}
\newcommand{\QQ}{\mathbb{Q}\tss}
\newcommand{\SSb}{\mathbb{S}\tss}
\newcommand{\ZZ}{\mathbb{Z}\tss}
\newcommand{\X}{ {\rm X}}
\newcommand{\Y}{ {\rm Y}}
\newcommand{\Z}{{\rm Z}}
\newcommand{\Ac}{\mathcal{A}}
\newcommand{\Lc}{\mathcal{L}}
\newcommand{\Mc}{\mathcal{M}}
\newcommand{\Pc}{\mathcal{P}}
\newcommand{\Qc}{\mathcal{Q}}
\newcommand{\Tc}{\mathcal{T}}
\newcommand{\Sc}{\mathcal{S}}
\newcommand{\Bc}{\mathcal{B}}
\newcommand{\Ec}{\mathcal{E}}
\newcommand{\Fc}{\mathcal{F}}
\newcommand{\Hc}{\mathcal{H}}
\newcommand{\Uc}{\mathcal{U}}
\newcommand{\Vc}{\mathcal{V}}
\newcommand{\Wc}{\mathcal{W}}
\newcommand{\Ar}{{\rm A}}
\newcommand{\Br}{{\rm B}}
\newcommand{\Ir}{{\rm I}}
\newcommand{\Fr}{{\rm F}}
\newcommand{\Jr}{{\rm J}}
\newcommand{\Or}{{\rm O}}
\newcommand{\Spr}{{\rm Sp}}
\newcommand{\Rr}{{\rm R}}
\newcommand{\Zr}{{\rm Z}}
\newcommand{\gl}{\mathfrak{gl}}
\newcommand{\middd}{{\rm mid}}
\newcommand{\Pf}{{\rm Pf}}
\newcommand{\Norm}{{\rm Norm\tss}}
\newcommand{\oa}{\mathfrak{o}}
\newcommand{\spa}{\mathfrak{sp}}
\newcommand{\osp}{\mathfrak{osp}}
\newcommand{\g}{\mathfrak{g}}
\newcommand{\h}{\mathfrak h}
\newcommand{\n}{\mathfrak n}
\newcommand{\z}{\mathfrak{z}}
\newcommand{\Zgot}{\mathfrak{Z}}
\newcommand{\p}{\mathfrak{p}}
\newcommand{\sll}{\mathfrak{sl}}
\newcommand{\agot}{\mathfrak{a}}
\newcommand{\qdet}{ {\rm qdet}\ts}
\newcommand{\Ber}{ {\rm Ber}\ts}
\newcommand{\HC}{ {\mathcal HC}}
\newcommand{\cdet}{ {\rm cdet}}
\newcommand{\tr}{ {\rm tr}}
\newcommand{\str}{ {\rm str}}
\newcommand{\loc}{{\rm loc}}
\newcommand{\Gr}{ {\rm Gr}\tss}
\newcommand{\sgn}{ {\rm sgn}\ts}
\newcommand{\ba}{\bar{a}}
\newcommand{\bb}{\bar{b}}
\newcommand{\bi}{\bar{\imath}}
\newcommand{\bj}{\bar{\jmath}}
\newcommand{\bk}{\bar{k}}
\newcommand{\bl}{\bar{l}}
\newcommand{\hb}{\mathbf{h}}
\newcommand{\Sym}{\mathfrak S}
\newcommand{\fand}{\quad\text{and}\quad}
\newcommand{\Fand}{\qquad\text{and}\qquad}

\renewcommand{\theequation}{\arabic{section}.\arabic{equation}}

\newtheorem{thm}{Theorem}[section]
\newtheorem{lem}[thm]{Lemma}
\newtheorem{prop}[thm]{Proposition}
\newtheorem{cor}[thm]{Corollary}
\newtheorem{conj}[thm]{Conjecture}
\newtheorem*{mthma}{Theorem A}
\newtheorem*{mthmb}{Theorem B}

\theoremstyle{definition}
\newtheorem{defin}[thm]{Definition}

\theoremstyle{remark}
\newtheorem{remark}[thm]{Remark}
\newtheorem{example}[thm]{Example}

\newcommand{\bth}{\begin{thm}}
\renewcommand{\eth}{\end{thm}}
\newcommand{\bpr}{\begin{prop}}
\newcommand{\epr}{\end{prop}}
\newcommand{\ble}{\begin{lem}}
\newcommand{\ele}{\end{lem}}
\newcommand{\bco}{\begin{cor}}
\newcommand{\eco}{\end{cor}}
\newcommand{\bde}{\begin{defin}}
\newcommand{\ede}{\end{defin}}
\newcommand{\bex}{\begin{example}}
\newcommand{\eex}{\end{example}}
\newcommand{\bre}{\begin{remark}}
\newcommand{\ere}{\end{remark}}
\newcommand{\bcj}{\begin{conj}}
\newcommand{\ecj}{\end{conj}}

\newcommand{\bal}{\begin{aligned}}
\newcommand{\eal}{\end{aligned}}
\newcommand{\beq}{\begin{equation}}
\newcommand{\eeq}{\end{equation}}
\newcommand{\ben}{\begin{equation*}}
\newcommand{\een}{\end{equation*}}

\newcommand{\bpf}{\begin{proof}}
\newcommand{\epf}{\end{proof}}

\def\beql#1{\begin{equation}\label{#1}}

\title{\Large\bf A new fusion procedure for the Brauer
algebra and evaluation homomorphisms}

\author{{A. P. Isaev,\quad A. I. Molev\quad and\quad
O. V. Ogievetsky\footnote{On leave of absence from
P. N. Lebedev Physical Institute, Leninsky Pr. 53,
117924 Moscow, Russia}}}

\date{} 
\maketitle


\begin{abstract}
We give a new fusion procedure for the Brauer algebra
by showing that all primitive idempotents
can be found by evaluating a rational
function in several variables which has the form
of a product of $R$-matrix type factors.
In particular, this provides a new fusion procedure
for the symmetric group involving an arbitrary parameter.
The $R$-matrices
are solutions of the Yang--Baxter equation
associated with the classical Lie algebras $\g_N$
of types $B$, $C$ and $D$.
Moreover, we construct an evaluation homomorphism
from a reflection equation algebra $\Br(\g_N)$ to $\U(\g_N)$
and show that the fusion procedure provides an equivalence
between natural tensor representations of $\Br(\g_N)$
with the corresponding evaluation modules.
\end{abstract}


\vspace{7 mm}

\noindent
Bogoliubov Laboratory of Theoretical Physics\newline
Joint Institute for Nuclear Research\newline
Dubna, Moscow region 141980, Russia\newline
isaevap@theor.jinr.ru

\vspace{7 mm}

\noindent
School of Mathematics and Statistics\newline
University of Sydney,
NSW 2006, Australia\newline
alexander.molev@sydney.edu.au

\vspace{7 mm}

\noindent
J.-V. Poncelet French--Russian Laboratory, UMI 2615 du CNRS\newline
Independent University of Moscow, 11 B. Vlasievski per.\newline
119002 Moscow, Russia\newline
{\small and}\newline
Center of Theoretical Physics\newline
Luminy, 13288 Marseille, France\newline
oleg@cpt.univ-mrs.fr

\newpage

\section{Introduction}
\label{sec:int}
\setcounter{equation}{0}

By an original observation of Jucys~\cite{j:yo},
all primitive idempotents
of the symmetric group $\Sym_n$ can be obtained by taking certain
limit values of the rational function
\beql{phiu}
\Phi(u_1,\dots,u_n)= \prod_{1\leqslant i<j\leqslant n}
\Big(1-\frac{s_{ij}}{u_i-u_j}\Big),
\eeq
where $s_{ij}\in\Sym_n$ is the transposition of $i$ and $j$,
$u_1,\dots,u_n$ are
complex variables and the product is calculated in the group algebra
$\CC[\Sym_n]$ in the lexicographic order on the pairs $(i,j)$.
This construction, which is commonly known as the {\it fusion
procedure\/}, was also developed by Cherednik~\cite{c:sb}, while
detailed proofs were given by Nazarov~\cite{n:yc}.
In the context of the quantum inverse scattering method
developed by Faddeev's Leningrad school,
the fusion procedure has been regarded as a way
to construct new solutions of the Yang--Baxter equation
out of old ones; see, e.g., \cite{krs:yb}.
The version of the fusion procedure found in \cite{m:fp} establishes
its equivalence to the construction of the idempotents
provided by Jucys~\cite{j:fy} and Murphy~\cite{m:nc}
in terms of some special elements of the group algebra of $\Sym_n$.
It was shown in \cite{imo:ih} that the fusion procedure for the Hecke
algebra admits a similar interpretation; cf. \cite{c:ni}, \cite{n:mh}.

A fusion procedure for the Brauer algebra $\Bc_n(\om)$ over $\CC(\om)$
was recently given by two of us in \cite{im:fp}.
For this algebra \eqref{phiu} was replaced by
the rational function
\beql{phiubra}
\Psi(u_1,\dots,u_n)=
\prod_{1\leqslant i<j\leqslant n}
\Big(1-\frac{\ep_{ij}}{u_i+u_j}\Big)
\prod_{1\leqslant i<j\leqslant n}
\Big(1-\frac{s_{ij}}{u_i-u_j}\Big)
\eeq
previously considered in \cite{n:rt},
with the ordered products as in \eqref{phiu};
the elements $\ep_{ij}$ and $s_{ij}$ of $\Bc_n(\om)$ are defined in
Sec.~\ref{sec:fpbra} below.

Recall that the irreducible
representations of $\Bc_n(\om)$
are indexed by all partitions of the nonnegative
integers $n,n-2,n-4,\dots$. If $\la$
is such a partition, then the
{\it updown $\la$-tableaux} $T$
parameterize basis vectors of the corresponding
representation; see Sec.~\ref{sec:fpbra} for the definitions.
This leads to an explicit isomorphism between $\Bc_n(\om)$
and the direct sum of matrix algebras. The {\it primitive
idempotents\/} $E^{\lambda}_{T}$ are the elements of $\Bc_n(\om)$
corresponding to the diagonal matrix units under this isomorphism;
see \cite{b:aw}, \cite{n:yo}, \cite{w:sb}.

By the main result of \cite{im:fp},
given an updown $\la$-tableau $T$,
the consecutive evaluations
\beql{reev}
(u_1-c_1)^{p_1}\dots (u_n-c_n)^{p_n}\ts
\Psi(u_1,\dots,u_n)\big|_{u_1=c_1}\big|_{u_2=c_2}\dots
\big|_{u_n=c_n}
\eeq
are well-defined and this value yields the
corresponding primitive idempotent
$E^{\lambda}_{T}$ multiplied by a nonzero constant
$f(T)$ which is calculated in an explicit form.
Here the $c_i$ are the contents of $T$
and $p_1,\dots,p_n$ are certain integers depending on $T$,
called its exponents.

The first main result of this paper
is a new construction
of all primitive idempotents of the Brauer algebra $\Bc_n(\om)$;
we use a different rational function in place of \eqref{phiubra}.
Namely, set
\begin{multline}\label{omubra}
\Omega(u_1,\dots,u_n)=
\prod_{1\leqslant i<j\leqslant n}
\Big(1+\frac{s_{ij}}{u_i+u_j-\om/2+1}-\frac{\ep_{ij}}{u_i+u_j}\Big)\\
{}\times\prod_{1\leqslant i<j\leqslant n}
\Big(1-\frac{s_{ij}}{u_i-u_j}+\frac{\ep_{ij}}{u_i-u_j-\om/2+1}\Big)
\end{multline}
with both products taken in the
lexicographic order on the pairs $(i,j)$. Thus,
\beql{omrm}
\Omega(u_1,\dots,u_n)=\prod_{1\leqslant i<j\leqslant n}
\rho_{ij}(-u_i-u_j+\vk)\ts \prod_{1\leqslant i<j\leqslant n}
\rho_{ij}(u_i-u_j),
\eeq
where we use the notation
\beql{rij}
\rho_{ij}(u)=1-\frac{s_{ij}}{u}+\frac{\ep_{ij}}{u-\vk},\qquad
\vk=\frac{\om}{2}-1.
\eeq
The rational functions $\rho_{ij}(u)$ satisfy
the Yang--Baxter equation
\beql{ybe}
\rho_{ij}(u)\ts \rho_{ik}(u+v)\ts \rho_{jk}(v)=
\rho_{jk}(v)\ts \rho_{ik}(u+v)\ts \rho_{ij}(u),
\eeq
with any distinct indices $i,j,k$, where we set
$\rho_{ji}(u)=\rho_{ij}(u)$ for $i<j$;
see \cite{zz:rf}. Also, for $i\ne j$ we have
\beql{unitarm}
\rho_{ij}(u)\ts\rho_{ij}(-u)=\frac{u^2-1}{u^2}.
\eeq
The new version of the fusion procedure
for $\Bc_n(\om)$ takes the following form: the consecutive
evaluations
\beql{reevom}
(u_1-c_1)^{p_1}\dots (u_n-c_n)^{p_n}\ts
\Omega(u_1,\dots,u_n)\big|_{u_1=c_1}\big|_{u_2=c_2}\dots
\big|_{u_n=c_n}
\eeq
are well-defined and this value yields the
product
$h(T)\tss E^{\lambda}_{T}$ as with \eqref{reev}, where $h(T)$
is a constant calculated in a way similar to $f(T)$.
To give an equivalent formulation, for any updown tableau
$T$ introduce a rational function $\Omega_T(u_1,\dots,u_n)$
which is obtained from $\Omega(u_1,\dots,u_n)$ by
multiplying by a numerical rational function in $u_1,\dots,u_n$
depending on the contents of $T$; see \eqref{omt}.
Then the fusion
procedure can be reformulated as the relation
\beql{reevomre}
E_T=\Omega_T(u_1,\dots,u_n)\big|_{u_1=c_1}\big|_{u_2=c_2}\dots
\big|_{u_n=c_n}.
\eeq

The existence of the two forms \eqref{reev} and \eqref{reevom}
(or, equivalently, \eqref{reevomre})
of the fusion procedure can be explained by their connections
with two different quantum algebras associated with a
classical Lie algebra of type $B$, $C$ or $D$.
Consider the natural action of the Brauer algebra $\Bc_n(\om)$ (with
an appropriate specialization
of the parameter $\om$)
on the tensor product space
\beql{cntenprint}
\CC^N\ot\CC^N\ot\dots\ot\CC^N,\qquad n\ \ \text{factors};
\eeq
see \eqref{braact} and
\eqref{braactsy} below. Using the expression for the idempotent
$E^{\lambda}_{T}$ provided by \eqref{reev} we find that
the subspace $E^{\lambda}_{T}(\CC^N)^{\ot n}$ carries
the structure of a representation of the {\it Olshanski twisted
Yangian\/} $\Y'(\g_N)$ associated with the orthogonal
Lie algebra $\g_N=\oa_N$
or symplectic Lie algebra $\g_N=\spa_N$, respectively,
($N$ is even for the latter) \cite{n:rt},
\cite{o:ty}; see also \cite[Ch.~2]{m:yc}.
On the other hand, the representation of
the universal enveloping algebra
$\U(\g_N)$ on the space $E^{\lambda}_{T}(\CC^N)^{\ot n}$
arising from the Brauer--Schur--Weyl duality extends to
the twisted Yangian $\Y'(\g_N)$ via the evaluation
homomorphism $\Y'(\g_N)\to\U(\g_N)$ \cite{o:ty}.
Thus, the fusion procedure in the form associated
with the evaluations \eqref{reev} yields an equivalence
of the two twisted Yangian actions on
the space $E^{\lambda}_{T}(\CC^N)^{\ot n}$.

In the new version of the fusion procedure associated
with the evaluations \eqref{reevom}, the role of the twisted
Yangian is now taken by the
{\it reflection algebra\/} $\Br(\g_N)$.
This algebra is defined by a reflection equation
and it is closely related to the
{\it Drinfeld Yangian\/} $\Y(\g_N)$; see Definition~\ref{def:refeq}
below.
The Yangian
$\Y(\g_N)$ contains the universal enveloping algebra
$\U(\g_N)$ as a subalgebra; however, due to
a result of Drinfeld~\cite{d:ha}, in contrast with the
Yangian for $\gl_N$,
there is no homomorphism $\Y(\g_N)\to \U(\g_N)$, identical
on the subalgebra $\U(\g_N)$.

The algebras $\Br(\oa_N)$ and $\Br(\spa_{2n})$ were
formally defined in \cite{aacfr:rp} together with their
super-version $\Br(\osp_{m|2n})$; however, they
do not seem to have received
much attention in the literature.
Our second main result (Theorem~\ref{thm:eval}) is a construction
of an evaluation homomorphism $\Br(\g_N)\to\U(\g_N)$.
Furthermore, we show that the algebra $\Br(\g_N)$ contains
the universal enveloping algebra $\U(\g_N)$ as a subalgebra,
and the evaluation homomorphism is identical on this
subalgebra. These results also apply to
the case of the Lie superalgebra $\osp_{m|2n}$, and we
indicate necessary changes in the notation in
Remark~\ref{rem:super}.

The expression for the idempotent
$E^{\lambda}_{T}$ provided by \eqref{reevom} indicates that
the subspace $E^{\lambda}_{T}(\CC^N)^{\ot n}$
carries the structure of a representation of the corresponding
reflection algebra $\Br(\g_N)$. Moreover, we show that,
as with the twisted Yangian, this representation
factors through the evaluation homomorphism $\Br(\g_N)\to\U(\g_N)$.

Taking the quotient of the Brauer algebra $\Bc_n(\om)$ by
the relations $\ep_{ij}=0$ we come to a new version of the fusion
procedure for the symmetric group $\Sym_n$ involving
an arbitrary parameter $\om$ (Corollary~\ref{cor:sym}).
The standard version of the procedure
associated with the rational function \eqref{phiu}
is recovered
in the limit $\om\to\infty$.
We show
that the new version is related to
the well-known evaluation homomorphism
$\Y(\gl_N)\to \U(\gl_N)$ and to its restriction to the subalgebra
of $\Y(\gl_N)$ isomorphic to a reflection algebra.

Both the fusion procedure and the evaluation homomorphism
described in Theorems~\ref{thm:newfus} and \ref{thm:eval}
admit their natural quantum analogues. Namely,
all primitive idempotents of the Birman--Murakami--Wenzl
algebra can be obtained by evaluating a universal rational
function. Moreover, it is possible to introduce
a $q$-analogue $\Br_q(\g_N)$ of the reflection algebra
and to construct a homomorphism from $\Br_q(\g_N)$
to the Drinfeld--Jimbo quantum group $\U_q(\g_N)$
with the properties similar to \eqref{evalya}.
The details will appear in our forthcoming publication
\cite{imo:fp}.

\medskip

A. P. I. and O. V. O. are grateful to the School of Mathematics and
Statistics of the University of Sydney for the warm hospitality
during their visits. We acknowledge the support of the Australian
Research Council. The work of A. P. I. was supported by the grants
RFBR 11-01-00980-a and RFBR-CNRS 07-02-92166-a.

\section{Fusion procedure}
\label{sec:fpbra}
\setcounter{equation}{0}

Let $n$ be a positive integer and $\om$ an indeterminate.
An $n$-diagram $d$ is a collection
of $2n$ dots arranged into two rows with $n$ dots in each row
connected by $n$ edges such that any dot belongs to only one edge.
The product of two diagrams $d_1$ and $d_2$ is determined by
placing $d_1$ above $d_2$ and identifying the vertices
of the bottom row of $d_1$ with the corresponding
vertices in the top row of $d_2$. Let $s$ be the number of
closed loops obtained in this placement. The product $d_1d_2$ is given by
$\om^{\tss s}$ times the resulting diagram without loops.
The {\it Brauer algebra\/} $\Bc_n(\om)$ is defined as the
$\CC(\om)$-linear span of the $n$-diagrams with the
multiplication defined above.
The dimension of the algebra is $1\cdot 3\cdots (2n-1)$.
The following presentation of $\Bc_n(\om)$ is well-known;
see, e.g., \cite{bw:bl}.

\bpr\label{prop:bradr}
The Brauer algebra $\Bc_n(\om)$ is isomorphic to the algebra
with $2n-2$ generators
$s_1,\dots,s_{n-1},\ep_1,\dots,\ep_{n-1}$
and the defining relations
\beql{bradr}
\begin{aligned}
s_i^2&=1,\qquad \ep_i^2=\om\ts \ee_i,\qquad
\su_i\ee_i=\ee_i\su_i=\ee_i,\qquad i=1,\dots,n-1,\\
s_i s_j &=s_js_i,\qquad \ee_i \ee_j = \ee_j \ee_i,\qquad
\su_i \ee_j = \ee_j\su_i,\qquad
|i-j|>1,\\
s_is_{i+1}s_i&=s_{i+1}s_is_{i+1},\qquad
\ee_i\ee_{i+1}\ee_i=\ee_i,\qquad \ee_{i+1}\ee_i\ee_{i+1}=\ee_{i+1},\\
\su_i\ee_{i+1}\ee_i&=\su_{i+1}\ee_i,\qquad \ee_{i+1}\ee_i
\su_{i+1}=\ee_{i+1}\su_i,\qquad
i=1,\dots,n-2.
\end{aligned}
\non
\end{equation}
\epr

The generators $s_i$ and $\ee_i$ correspond to the following diagrams
respectively:

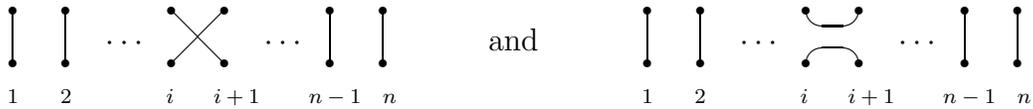
\begin{center}
\begin{picture}(400,60)
\thinlines

\put(10,20){\circle*{3}}
\put(30,20){\circle*{3}}
\put(70,20){\circle*{3}}
\put(90,20){\circle*{3}}
\put(130,20){\circle*{3}}
\put(150,20){\circle*{3}}

\put(10,40){\circle*{3}}
\put(30,40){\circle*{3}}
\put(70,40){\circle*{3}}
\put(90,40){\circle*{3}}
\put(130,40){\circle*{3}}
\put(150,40){\circle*{3}}

\put(10,20){\line(0,1){20}}
\put(30,20){\line(0,1){20}}
\put(70,20){\line(1,1){20}}
\put(90,20){\line(-1,1){20}}
\put(130,20){\line(0,1){20}}
\put(150,20){\line(0,1){20}}

\put(45,25){$\cdots$}
\put(105,25){$\cdots$}

\put(8,5){\scriptsize $1$ }
\put(28,5){\scriptsize $2$ }
\put(68,5){\scriptsize $i$ }
\put(86,5){\scriptsize $i+1$ }
\put(122,5){\scriptsize $n-1$ }
\put(150,5){\scriptsize $n$ }

\put(190,25){\text{and}}

\put(250,20){\circle*{3}}
\put(270,20){\circle*{3}}
\put(310,20){\circle*{3}}
\put(330,20){\circle*{3}}
\put(370,20){\circle*{3}}
\put(390,20){\circle*{3}}

\put(250,40){\circle*{3}}
\put(270,40){\circle*{3}}
\put(310,40){\circle*{3}}
\put(330,40){\circle*{3}}
\put(370,40){\circle*{3}}
\put(390,40){\circle*{3}}

\put(250,20){\line(0,1){20}}
\put(270,20){\line(0,1){20}}
\put(320,20){\oval(20,12)[t]}
\put(320,40){\oval(20,12)[b]}
\put(370,20){\line(0,1){20}}
\put(390,20){\line(0,1){20}}

\put(285,25){$\cdots$}
\put(345,25){$\cdots$}

\put(248,5){\scriptsize $1$ }
\put(268,5){\scriptsize $2$ }
\put(308,5){\scriptsize $i$ }
\put(326,5){\scriptsize $i+1$ }
\put(362,5){\scriptsize $n-1$ }
\put(390,5){\scriptsize $n$ }

\end{picture}
\end{center}

The subalgebra of $\Bc_n(\om)$ generated over $\CC$
by $s_1,\dots,s_{n-1}$
is isomorphic to the group algebra $\CC[\Sym_n]$ so that $s_i$
can be identified with the transposition $(i,i+1)$.
Then for any $1\leqslant i<j\leqslant n$
the transposition $s_{ij}=(i,j)$ can be regarded
as an element of $\Bc_n(\om)$. Moreover, $\ep_{ij}$
will denote the element of $\Bc_n(\om)$ represented
by the diagram in which the $i$-th and $j$-th dots in the top row,
as well as the $i$-th and $j$-th dots in the bottom
row are connected by an edge, while
the remaining edges
connect the $k$-th dot in the top row
with the $k$-th dot in the bottom row for each $k\ne i,j$.
Equivalently, in terms of the presentation of $\Bc_n(\om)$
provided by Proposition~\ref{prop:bradr},
\ben
s_{ij}=s_i\tss s_{i+1}\dots s_{j-2}\tss
s_{j-1}\tss s_{j-2}\dots s_{i+1}\tss s_i
\Fand
\ep_{ij}=s_{i,j-1}\tss \ep_{j-1}\tss s_{i,j-1}.
\een
We also set $\ep_{ji}=\ep_{ij}$ and $s_{ji}=s_{ij}$ for $i<j$.
The Brauer algebra $\Bc_{n-1}(\om)$ can be regarded
as the subalgebra of $\Bc_n(\om)$ spanned by all diagrams
in which the $n$-th dots in the top and bottom rows are
connected by an edge.

The {\it Jucys--Murphy elements\/} $x_1,\dots,x_n$
for the Brauer algebra $\Bc_n(\om)$ are given by the formulas
\beql{jmdef}
x_r=\frac{\om-1}{2}+\sum_{i=1}^{r-1}(s_{i\tss r}-\ep_{i\tss r}),
\qquad r=1,\dots,n;
\eeq
see \cite{lr:rh} and \cite{n:yo}, where, in particular,
the eigenvalues
of the $x_r$ in irreducible representations were calculated.
The element $x_n$ commutes
with the subalgebra $\Bc_{n-1}(\om)$. This implies that
the elements $x_1,\dots,x_n$
of $\Bc_n(\om)$ pairwise commute. They can be used
to construct a complete set of pairwise orthogonal primitive
idempotents for the Brauer algebra following
the approach of Jucys~\cite{j:fy} and Murphy~\cite{m:nc}.
Namely, let $\la$ be a partition of $n-2f$
for some $f\in\{0,1,\dots,\lfloor n/2\rfloor\}$.
We will identify partitions with their diagrams
so that if the parts of $\la$ are $\la_1,\la_2,\dots$ then
the corresponding diagram is
a left-justified array of rows of unit boxes containing
$\la_1$ boxes in the top row, $\la_2$ boxes
in the second row, etc.
The box in row $i$ and column $j$ of a diagram
will be denoted as the pair $(i,j)$.
An {\it updown $\la$-tableau\/} is a sequence
$T=(\La_1,\dots,\La_n)$ of diagrams such that
for each $r=1,\dots,n$ the diagram $\La_r$ is obtained
from $\La_{r-1}$ by adding or removing one box,
where we set $\La_0=\varnothing$, the empty diagram, and $\La_n=\la$.
To each updown tableau $T$ we attach the corresponding
sequence of {\it contents\/} $(c_1,\dots,c_n)$, $c_r=c_r(T)$,
where
\ben
c_r=\frac{\om-1}{2}+j-i\qquad\text{or}\qquad
c_r=-\Big(\frac{\om-1}{2}+j-i\Big),
\een
if $\La_r$ is obtained by adding the box $(i,j)$ to $\La_{r-1}$
or by removing this box from $\La_{r-1}$,
respectively. The primitive idempotents $E_{T}=E^{\lambda}_{T}$
can now be defined by
the following recurrence formula (we omit the superscripts
indicating the diagrams since they are determined by
the updown tableaux).
Set $\mu=\La_{n-1}$ and consider the updown $\mu$-tableau
$U=(\La_1,\dots,\La_{n-1})$. Let $\alpha$ be the box
which is added to or removed from $\mu$ to get $\la$. Then
\beql{murphyfo}
E_{T}=E_{U}\ts
\frac{(x_n-a_1)\dots (x_n-a_k)}{(c_n-a_1)\dots (c_n-a_k)},
\eeq
where $a_1,\dots,a_k$ are the contents of all boxes
excluding $\alpha$,
which can be removed from or added
to $\mu$ to get a diagram.
When $\la$ runs over all partitions
of $n,n-2,\dots$, and $T$ runs over all updown $\la$-tableaux,
the elements $\{E_{T}\}$ yield a complete set
of pairwise orthogonal primitive
idempotents for $\Bc_n(\om)$. They have the properties
\beql{xiet}
x_r\ts E_{T}=E_{T}\ts x_r=c_r(T)\ts E_{T}, \qquad r=1,\dots,n.
\eeq
Moreover, given an updown tableau $U=(\La_1,\dots,\La_{n-1})$,
we have the relation
\beql{eutet}
E_{U}=\sum_{T} E_{T},
\eeq
the sum is over all updown tableaux
of the form $T=(\La_1,\dots,\La_{n-1},\La_n)$; see
e.g. \cite{lr:rh} and \cite{n:yo}.
Relation \eqref{murphyfo} implies
\beql{jmform}
E_{T}=E_{U}\ts \frac{(u-c_n)(u+x_n-\vk)}
{(u-x_n)(u+c_n-\vk)}\ts
\Big|^{}_{u=c_n},
\eeq
where $u$ is a complex variable and we use notation \eqref{rij}.
This relation follows
by application of \eqref{eutet} and then \eqref{xiet}.

Given an updown tableau $T=(\La_1,\dots,\La_n)$
with the respective contents $c_1,\dots,c_n$,
introduce a rational function in $u_1,\dots,u_n$
with values in the Brauer algebra $\Bc_n(\om)$ by
\beql{omt}
\Omega_T(u_1,\dots,u_n)=
\prod_{r=2}^{n}\frac{(u_r-c_r)(u_r+c_1-\vk)}{(u_r-c_1)(u_r+c_r-\vk)}
\ts\prod_{i=1}^{r-1}\frac{(u_r-u_i)^2}{(u_r-u_i)^2-1}\ts
\Omega(u_1,\dots,u_n),
\eeq
where $\Omega(u_1,\dots,u_n)$ is defined in \eqref{omubra}
and \eqref{omrm}. Note that if
the indices $i,j,k,l$ are distinct then
the elements $\rho_{ij}(u)$ and $\rho_{kl}(v)$ commute.
Together with \eqref{ybe} this implies
the relation
\begin{multline}\non
\rho_{1,n}(-u_1-u_n+\vk)\dots \rho_{n-1,n}(-u_{n-1}-u_n+\vk)\ts
\prod_{1\leqslant i<j\leqslant n-1}
\rho_{ij}(u_i-u_j)\\
{}=\prod_{1\leqslant i<j\leqslant n-1}
\rho_{ij}(u_i-u_j)\ts \rho_{n-1,n}(-u_{n-1}-u_n+\vk)\dots
\rho_{1,n}(-u_1-u_n+\vk).
\end{multline}
An easy induction on $n$ leads to the following
equivalent expression for the rational function
$\Omega(u_1,\dots,u_n)$:
\begin{multline}\label{psedec}
\Omega(u_1,\dots,u_n)=\prod_{r=2}^{n}\ts
\rho_{r-1,r}(-u_{r-1}-u_r+\vk)\dots \rho_{1,r}(-u_1-u_r+\vk)\\
{}\times \rho_{1,r}(u_1-u_r)
\dots \rho_{r-1,r}(u_{r-1}-u_r),
\end{multline}
where the factors are ordered in accordance with
the increasing values of $r$.

\bth\label{thm:newfus}
The idempotent $E_T$ is found by
the consecutive evaluations
\ben
E_T=\Omega_T(u_1,\dots,u_n)\big|_{u_1=c_1}\big|_{u_2=c_2}\dots
\big|_{u_n=c_n}.
\een
\eth

\bpf
We use the induction on $n$.
By the induction hypothesis, setting $u=u_n$
and using \eqref{psedec} we get
\begin{multline}\label{indugeth}
\Omega_T(u_1,\dots,u_n)\big|_{u_1=c_1}\big|_{u_2=c_2}\dots
\big|_{u_{n-1}=c_{n-1}}
=\frac{(u-c_n)(u+c_1-\vk)}{(u-c_1)(u+c_n-\vk)}\ts
\ts\prod_{i=1}^{n-1}\frac{(u-c_i)^2}{(u-c_i)^2-1}\ts\\
{}\times{}E_U\ts
\rho_{n-1,n}(-c_{n-1}-u+\vk)\dots \rho_{1,n}(-c_1-u+\vk)\ts
\rho_{1,n}(c_1-u)\dots \rho_{n-1,n}(c_{n-1}-u),
\end{multline}
where $U$ is the updown tableau $(\La_1,\dots,\La_{n-1})$.

\ble\label{lem:jmsi}
We have the identity
\begin{multline}\label{jm}
E_U\ts
\rho_{n-1,n}(-c_{n-1}-u+\vk)\dots \rho_{1,n}(-c_1-u+\vk)\ts
\rho_{1,n}(c_1-u)\dots \rho_{n-1,n}(c_{n-1}-u)\\
{}=\frac{u-c_1}{u+c_1-\vk}\ts\prod_{i=1}^{n-1}
\ts\frac{(u-c_i)^2-1}{(u-c_i)^2}
\ts E_U\ts\frac{u+x_n-\vk}{u-x_n}.
\end{multline}
\ele

\bpf
Embed the Brauer algebra $\Bc_n(\om)$ into $\Bc_m(\om)$
for some $m\geqslant n$ and
prove a more general identity
\begin{multline}\label{jmm}
E_U\ts
\rho_{n-1,m}(-c_{n-1}-u+\vk)\dots \rho_{1,m}(-c_1-u+\vk)\ts
\rho_{1,m}(c_1-u)\dots \rho_{n-1,m}(c_{n-1}-u)\\
{}=\frac{u-c_1}{u+c_1-\vk}\ts\prod_{i=1}^{n-1}
\ts\frac{(u-c_i)^2-1}{(u-c_i)^2}
\ts E_U\ts\frac{u+x^{(n-1)}_m-\vk}{u-x^{(n-1)}_m},
\end{multline}
where
\beql{jmdefm}
x^{(k)}_m=\frac{\om-1}{2}+\sum_{i=1}^{k}(s_{i\tss m}-\ep_{i\tss m});
\eeq
in particular, $x^{(n-1)}_n=x_n$; see \eqref{jmdef}. We use induction
on $n$ while keeping $m$ fixed. In the case $n=1$ the identity
is certainly true.
Suppose that $n\geqslant 2$. By \eqref{eutet} we
have $E_U=E_U\ts E_W$,
where $W$ is the updown tableau $(\La_1,\dots,\La_{n-2})$
(in the case $n=2$ we set $E_W=1$).
Hence, using the induction hypothesis
we can write the left hand side of \eqref{jmm} as
\begin{multline}
E_U\ts
\rho_{n-1,m}(-c_{n-1}-u+\vk)\\
{}\times
\frac{u-c_1}{u+c_1-\vk}\ts\prod_{i=1}^{n-2}
\ts\frac{(u-c_i)^2-1}{(u-c_i)^2}
\ts E_W\ts\frac{u+x^{(n-2)}_m-\vk}{u-x^{(n-2)}_m}
\ts \rho_{n-1,m}(c_{n-1}-u)
\non
\end{multline}
which equals
\beql{lhstra}
\frac{u-c_1}{u+c_1-\vk}\ts\prod_{i=1}^{n-2}
\ts\frac{(u-c_i)^2-1}{(u-c_i)^2}\ts E_U\ts
\rho_{n-1,m}(-c_{n-1}-u+\vk)\ts
\frac{u+x^{(n-2)}_m-\vk}{u-x^{(n-2)}_m}
\ts \rho_{n-1,m}(c_{n-1}-u).
\eeq
Applying
\eqref{unitarm} we find that
\ben
\rho_{n-1,m}(c_{n-1}-u)=\frac{(u-c_{n-1})^2-1}{(u-c_{n-1})^2}\ts
\rho_{n-1,m}(u-c_{n-1})^{-1}.
\een
Hence, comparing \eqref{lhstra} with the right hand side of
\eqref{jmm}, we conclude that the lemma will be implied
by the identity
\ben
E_U\ts
\rho_{n-1,m}(-c_{n-1}-u+\vk)\ts\frac{u+x^{(n-2)}_m-\vk}{u-x^{(n-2)}_m}
=E_U\ts \frac{u+x^{(n-1)}_m-\vk}{u-x^{(n-1)}_m}\ts
\rho_{n-1,m}(u-c_{n-1}).
\een
Since $x^{(n-1)}_m$ commutes with $E_U$, we can write the
identity in the equivalent form
\begin{multline}
E_U\ts(u-x^{(n-1)}_m)\ts\Big(1+\frac{s_{n-1,m}}{c_{n-1}+u-\vk}-
\frac{\ep_{n-1,m}}{c_{n-1}+u}\Big)\ts (u+x^{(n-2)}_m-\vk)\\
{}-E_U\ts(u+x^{(n-1)}_m-\vk)\ts\Big(1+\frac{s_{n-1,m}}{c_{n-1}-u}-
\frac{\ep_{n-1,m}}{c_{n-1}-u+\vk}\Big)\ts (u-x^{(n-2)}_m)=0.
\non
\end{multline}
To verify the latter, note that by \eqref{xiet} and
the relations
in the Brauer algebra,
\begin{multline}
E_U\ts x^{(n-1)}_m \ts s_{n-1,m}\ts x^{(n-2)}_m
=E_U\ts x^{(n-1)}_m \ts x_{n-1}\ts s_{n-1,m}\\
=E_U\ts x_{n-1}\ts x^{(n-1)}_m \ts s_{n-1,m}
=c_{n-1}\ts E_U\ts x^{(n-1)}_m \ts s_{n-1,m},
\non
\end{multline}
while
\ben
E_U\ts x^{(n-1)}_m \ts \ep_{n-1,m}\ts x^{(n-2)}_m
=-E_U\ts x_{n-1}\ts \ep_{n-1,m}\ts x^{(n-2)}_m
=-c_{n-1}\ts E_U\ts \ep_{n-1,m}\ts x^{(n-2)}_m.
\een
Together with the relation $x^{(n-1)}_m=x^{(n-2)}_m+s_{n-1,m}
-\ep_{n-1,m}$ this yields the identity in question and
completes the proof of the lemma.
\epf

Return to the proof of the theorem; we are left to
verify
that the rational function
\ben
E_{U}\ts \frac{(u-c_n)(u+x_n-\vk)}
{(u-x_n)(u+c_n-\vk)}
\een
is well-defined at $u=c_n$ and the value coincides with $E_T$.
But this holds due to \eqref{jmform}.
\epf

\bex\label{exam:idem}\quad (i)\ \
In the case of the Brauer algebra $\Bc_2(\om)$ we have
three updown tableaux

\setlength{\unitlength}{0.3em}
\begin{center}
\begin{picture}(82,5)

\put(-10,0){$U_1=\big($}

\put(0,0){\line(0,1){2}}
\put(2,0){\line(0,1){2}}
\put(0,0){\line(1,0){2}}
\put(0,2){\line(1,0){2}}
\put(2.5,0){,}

\put(8,0){\line(0,1){2}}
\put(10,0){\line(0,1){2}}
\put(12,0){\line(0,1){2}}
\put(8,0){\line(1,0){4}}
\put(8,2){\line(1,0){4}}
\put(12.5,0)

\put(13,0){$\big),$}

\put(26,0){$U_2=\big($}

\put(36,0){\line(0,1){2}}
\put(38,0){\line(0,1){2}}
\put(36,0){\line(1,0){2}}
\put(36,2){\line(1,0){2}}
\put(38.5,0){,}

\put(44,-2){\line(0,1){4}}
\put(46,-2){\line(0,1){4}}
\put(44,-2){\line(1,0){2}}
\put(44,0){\line(1,0){2}}
\put(44,2){\line(1,0){2}}
\put(46.5,0)

\put(48,0){$\big),$}

\put(61,0){$U_3=\big($}

\put(71,0){\line(0,1){2}}
\put(73,0){\line(0,1){2}}
\put(71,0){\line(1,0){2}}
\put(71,2){\line(1,0){2}}
\put(73.5,0){,}

\put(79,0){$\varnothing$}

\put(82,0){$\big).$}

\end{picture}
\end{center}
\setlength{\unitlength}{1pt}

\noindent
The respective contents $(c_1,c_2)$ are given by
\ben
\Big(\frac{\om-1}{2},\ts\frac{\om+1}{2}\Big),\qquad
\Big(\frac{\om-1}{2},\ts\frac{\om-3}{2}\Big),\qquad
\Big(\frac{\om-1}{2},\ts-\frac{\om-1}{2}\Big).
\een
The definition \eqref{omt} reads
\ben
\bal
\Omega_T(u_1,u_2)&=\frac{(u_2-c_2)(u_2+c_1-\vk)}{(u_2-c_1)(u_2+c_2-\vk)}
\ts\frac{(u_1-u_2)^2}{(u_1-u_2)^2-1}\ts\\
{}&\times\Big(1+\frac{s_1}{u_1+u_2-\vk}-\frac{\ep_1}{u_1+u_2}\Big)
\ts \Big(1-\frac{s_1}{u_1-u_2}+\frac{\ep_1}{u_1-u_2-\vk}\Big)
\eal
\een
so that applying Theorem~\ref{thm:newfus}, we find
\ben
E_{U_1}=
\frac{1+s_1}{2}-\frac{\ep_1}{\om},\qquad
E_{U_2}=\frac{1-s_1}{2},\qquad
E_{U_3}=\frac{\ep_1}{\om}.
\een

\quad (ii)\ \
In the case of $\Bc_3(\om)$ consider the following
two updown tableaux
associated with the diagram $\la=(1)$:
\setlength{\unitlength}{0.3em}
\begin{center}
\begin{picture}(65,5)

\put(-10,0){$T_1=\big($}

\put(0,0){\line(0,1){2}}
\put(2,0){\line(0,1){2}}
\put(0,0){\line(1,0){2}}
\put(0,2){\line(1,0){2}}
\put(2.5,0){,}

\put(8,-2){\line(0,1){4}}
\put(10,-2){\line(0,1){4}}
\put(8,-2){\line(1,0){2}}
\put(8,0){\line(1,0){2}}
\put(8,2){\line(1,0){2}}
\put(10.5,0){,}

\put(17,0){\line(0,1){2}}
\put(19,0){\line(0,1){2}}
\put(17,0){\line(1,0){2}}
\put(17,2){\line(1,0){2}}

\put(20,0){$\big),$}

\put(33,0){$T_2=\big($}

\put(43,0){\line(0,1){2}}
\put(45,0){\line(0,1){2}}
\put(43,0){\line(1,0){2}}
\put(43,2){\line(1,0){2}}
\put(45.5,0){,}

\put(51,0){$\varnothing,$}

\put(60,0){\line(0,1){2}}
\put(62,0){\line(0,1){2}}
\put(60,0){\line(1,0){2}}
\put(60,2){\line(1,0){2}}

\put(62.5,0){$\big).$}

\end{picture}
\end{center}
\setlength{\unitlength}{1pt}

\noindent
The respective contents $(c_1,c_2,c_3)$ are given by
\ben
\Big(\frac{\om-1}{2},\ts\frac{\om-3}{2},\ts-\frac{\om-3}{2}\Big),\qquad
\Big(\frac{\om-1}{2},\ts-\frac{\om-1}{2},
\ts\frac{\om-1}{2}\Big).
\een
Then using Theorem~\ref{thm:newfus} and \eqref{indugeth} we find
that $E_{T_1}$ and $E_{T_2}$
can be calculated by evaluating the respective rational
functions
\begin{multline}
\frac{(u-c_3)(u+c_1-\vk)}{(u-c_1)(u+c_3-\vk)}\ts
\ts\frac{(u-c_1)^2}{(u-c_1)^2-1}\ts\frac{(u-c_2)^2}{(u-c_2)^2-1}\ts\\
{}\times{}E_{U}\ts
\rho^{}_{23}(-c_2-u+\vk)\ts\rho^{}_{13}(-c_1-u+\vk)\ts
\rho^{}_{13}(c_1-u)\ts\rho^{}_{23}(c_2-u)
\non
\end{multline}
at $u=c_3$ with $U=U_2$ and $U=U_3$, respectively;
see Example~(i).
Performing
the calculation we get
\ben
E_{T_1}=
\frac{(1-s_1)\ts \ep_2\ts(1-s_1)}{2\ts(\om-1)},\qquad
E_{T_2}=\frac{\ep_1}{\om}.
\een

\quad (iii)\ \
For $\Bc_n(\om)$ consider the only updown tableau
$T=\big((1),(2),\dots,(n)\big)$
associated with the diagram $\la=(n)$.
The corresponding contents $(c_1,\dots,c_n)$ are
found by the formula $c_i=(\om-1)/2+i-1$.
Hence, by Theorem~\ref{thm:newfus}
the idempotent $S_n=E_T$
(the {\it symmetrizer\/} in the Brauer algebra) is given by
\ben
S_n=\frac{1}{n!}\ts\prod_{r=2}^n\ts\frac{\om+2r-2}{\om+4r-4}
\prod_{1\leqslant i<j\leqslant n}
\rho_{ij}(2-i-j-\om/2)\ts \prod_{1\leqslant i<j\leqslant n}
\rho_{ij}(i-j)
\een
with the products over the pairs $(i,j)$ taken in the
lexicographic order.
Given any index $k\in\{1,\dots,n-1\}$,
we can use \eqref{ybe} to reorder the factors in
the last product in such a way
that it would begin with the factor $\rho_{k,k+1}(-1)$.
Since
\ben
\ep_k\ts \rho_{k,k+1}(-1)=0\Fand
s_k\ts \rho_{k,k+1}(-1)=\rho_{k,k+1}(-1),
\een
we find that for any $k<l$
\ben
\rho_{kl}(u)\ts \prod_{1\leqslant i<j\leqslant n}
\rho_{ij}(i-j)=\frac{u-1}{u}\ts\prod_{1\leqslant i<j\leqslant n}
\rho_{ij}(i-j).
\een
Hence the expression for the symmetrizer simplifies to
\ben
S_n=\frac{1}{n!}\ts \prod_{1\leqslant i<j\leqslant n}
\rho_{ij}(i-j);
\een
cf. \cite[Remark~3.8]{im:fp}.

\quad (iv)\ \
In the case of $\Bc_n(\om)$ consider the only updown tableau
$T=\big((1),(1^2),\dots,(1^n)\big)$
associated with the diagram $\la=(1^n)$.
The corresponding contents $(c_1,\dots,c_n)$ are
given by $c_i=(\om-1)/2-i+1$. Hence, by Theorem~\ref{thm:newfus}
the idempotent $A_n=E_T$
(the {\it anti-symmetrizer\/} in the Brauer algebra) is given by
\beql{newfp}
A_n=\frac{1}{n!}\ts\prod_{r=2}^n\ts\frac{\om-2r+2}{\om-4r+4}
\prod_{1\leqslant i<j\leqslant n}
\rho_{ij}(i+j-2-\om/2)\ts \prod_{1\leqslant i<j\leqslant n}
\rho_{ij}(j-i)
\eeq
with the products over the pairs $(i,j)$ taken in the
lexicographic order.
A different expression
for the anti-symmetrizer is provided by the alternative
fusion procedure in \cite[Remark~3.8]{im:fp}:
\beql{oldfp}
A_n=\frac{1}{n!}\ts\prod_{1\leqslant i<j\leqslant n}
\Big(1-\frac{s_{ij}}{j-i}\Big)
\eeq
with the product taken in the
lexicographic order on the pairs $(i,j)$. A more direct way to
see the coincidence of the elements given by \eqref{newfp}
and \eqref{oldfp} is to regard the group algebra $\CC[\Sym_n]$
for the symmetric group as a natural subalgebra of $\Bc_n(\om)$
and observe that \eqref{oldfp} is the well-known factorization
of the anti-symmetrizer in $\CC[\Sym_n]$. Indeed,
these anti-symmetrizers satisfy the well-known recurrence
relation $n\ts A_n=A_{n-1}\big(1-s_{1n}-\dots-s_{n-1,n}\big)$.
Therefore, employing the induction
on $n$ and using the recurrence relations \eqref{murphyfo}
in $\Bc_n(\om)$
we come to checking the identity
\ben
A_{n-1}\frac{(x_n-c_1-1)(x_n+c_1-n+2)}{2\tss c_1-2n+3}
=A_{n-1}\big({-1}+s_{1n}+\dots+s_{n-1,n}\big),
\een
where, as before, $c_1=(\om-1)/2$. The calculation
is straightforward and it relies on the
following relations
in $\Bc_n(\om)$ which hold for any two distinct
indices $i,j\in\{1,\dots,n-1\}$:
\ben
A_{n-1}\ts s_{in}\ts s_{jn}=-A_{n-1}\ts s_{in},\quad
A_{n-1}\ts s_{in}\ts \ep_{jn}=0,\quad
A_{n-1}\ts \ep_{in}\ts \ep_{jn}=A_{n-1}\ts \ep_{in}\ts s_{ij}=
-A_{n-1}\ts \ep_{jn}
\een
and
\ben
A_{n-1}(\ep_{in} s_{jn} + \ep_{jn} s_{in}) =
A_{n-1}(\ep_{in} s_{jn} - s_{ij} \ep_{jn} s_{in}) =
A_{n-1}(\ep_{in} s_{jn} - \ep_{in} s_{jn}) = 0.
\een
\vskip-1.6\baselineskip
\qed
\eex

To make a connection of the version of the fusion procedure
provided by Theorem~\ref{thm:newfus} with that of \cite{im:fp},
recall some parameters associated with
updown tableaux; see \cite{im:fp}.
Given an updown $\mu$-tableau
$U=(\La_1,\dots,\La_{n-1})$ we define two infinite matrices
$m(U)$ and $m'(U)$ whose rows and columns are labelled
by positive integers and only a finite number of entries
in each of the matrices are nonzero. The entry $m_{ij}$
of the matrix $m(U)$ (resp., the entry $m'_{ij}$
of the matrix $m'(U)$) equals the number of times
the box $(i,j)$ was added (resp., removed) in the sequence
of diagrams $(\varnothing=\La_0,\La_1,\dots,\La_{n-1})$.
For each integer $k$
define the nonnegative integers $d_k=d_k(U)$ and
$d^{\tss\prime}_k=d^{\tss\prime}_k(U)$
as the respective sums of the entries of the matrices
$m(U)$ and $m'(U)$ on the $k$-th diagonal,
\ben
d_k=\sum_{j-i=k}m_{ij},\qquad d^{\tss\prime}_k=\sum_{j-i=k}m'_{ij},
\een
and set
\beql{gk}
g_k(U)=\de_{k\tss 0}+d_{k-1}+d_{k+1}-2\tss d_k,\qquad
g'_k(U)=d^{\tss\prime}_{k-1}
+d^{\tss\prime}_{k+1}-2\tss d^{\tss\prime}_k.
\eeq
Now the {\it exponents\/} $p_1,\dots,p_n$ of
an updown $\la$-tableau $T=(\La_1,\dots,\La_n)$ are defined
inductively, so that $p_r$ depends only on the first $r$
diagrams $(\La_1,\dots,\La_r)$ of $T$. Hence, it is sufficient
to define $p_n$. Taking $U=(\La_1,\dots,\La_{n-1})$ we set
\beql{pn}
p_n=1-g_{k_n}(U)\qquad\text{or}\qquad p_n=1-g'_{k_n}(U),
\eeq
respectively, if $\La_n$ is obtained from $\La_{n-1}$
by adding a box on the diagonal $k_n$ or by
removing a box on the diagonal $k_n$.

Define the constants $h(T)$ inductively by the formula
\beql{ft}
h(T)=h(U)\ts\psi(U,T),
\eeq
where $U=(\La_1,\dots,\La_{n-1})$,
$T=(\La_1,\dots,\La_n)$, and $h(U)=1$ if $U=(\La_1)$.
Here
\ben
\psi(U,T)=\frac{4k_n+\om}{2k_n+\om}\ts\prod_{k\ne k_n}(k_n-k)^{g_k}
\prod_k(k_n+k+\om-1)^{g^{\tss\prime}_k}
\een
or
\ben
\psi(U,T)=\frac{4k_n+6\om-4}{2k_n+\om-2}\ts
\prod_{k\ne k_n}(-k_n+k)^{g^{\tss\prime}_k}
\prod_k(-k_n-k-\om+1)^{g_k},
\een
if $\La_n$ is obtained from $\La_{n-1}$
by adding or removing a box on the diagonal $k_n$,
respectively, where the products are taken
over all integers $k$, while $g_k=g_k(U)$ and $g'_k=g'_k(U)$.

Consider now the rational function $\Omega(u_1,\dots,u_n)$
with values in the Brauer algebra $\Bc_n(\om)$
defined by \eqref{omubra}.

\bco\label{cor:newfuseq}
For any updown tableau $T=(\La_1,\dots,\La_n)$
the consecutive evaluations
\ben
(u_1-c_1)^{p_1}\dots (u_n-c_n)^{p_n}\ts
\Omega(u_1,\dots,u_n)\big|_{u_1=c_1}\big|_{u_2=c_2}\dots
\big|_{u_n=c_n}
\een
are well-defined. The corresponding value coincides with
$h(T)\tss E_{T}$.
\eco

\bpf
We use induction on $n$ as in the proof of Theorem~\ref{thm:newfus}.
By the induction hypothesis, setting $u=u_n$ we get
\begin{multline}\label{induge}
(u_1-c_1)^{p_1}\dots (u_n-c_n)^{p_n}\ts
\Omega(u_1,\dots,u_n)\big|_{u_1=c_1}\big|_{u_2=c_2}\dots
\big|_{u_{n-1}=c_{n-1}}=h(U)\tss (u-c_n)^{p_n}\\
{}\times{}E_U\ts
\rho_{n-1,n}(-c_{n-1}-u+\vk)\dots \rho_{1,n}(-c_1-u+\vk)\ts
\rho_{1,n}(c_1-u)\dots \rho_{n-1,n}(c_{n-1}-u),
\end{multline}
where $U$ is the updown tableau $(\La_1,\dots,\La_{n-1})$.
Applying Lemma~\ref{lem:jmsi}, we come to
showing that the rational function
\ben
h(U)\tss
\frac{(u-c_1)(u+c_n-\vk)}{u+c_1-\vk}\ts
\prod_{r=1}^{n-1}\Big(1-\frac{1}{(u-c_r)^2}\Big)
\ts(u-c_n)^{p_n-1}\cdot
E_U\ts \frac{u-c_n}{u-x_n}\ts \frac{u+x_n-\vk}{u+c_n-\vk}
\een
is regular at $u=c_n$ and its value equals $h(T)\tss E_T$.
Using the parameters \eqref{gk}, we can write
this expression as
\begin{multline}
h(U)\ts\frac{u+c_n-\vk}{u+c_1-\vk}\ts
\prod_{k}\Big(u-\frac{\om-1}{2}-k\Big)^{g_k}
\prod_{k}\Big(u+\frac{\om-1}{2}+k\Big)^{g^{\tss\prime}_k}\ts
(u-c_n)^{p_n-1}\\
{}\times E_U\ts \frac{u-c_n}{u-x_n}\ts \frac{u+x_n-\vk}{u+c_n-\vk},
\non
\end{multline}
where $k$ runs over the set of integers.
If the diagram $\La_n$ is obtained from $\La_{n-1}$
by adding or removing a box on the diagonal $k_n$, then
the value of the content $c_n$ is given by the respective formulas
\ben
c_n=\frac{\om-1}{2}+k_n\qquad\text{or}\qquad
c_n=-\Big(\frac{\om-1}{2}+k_n\Big).
\een
The argument is completed by recalling the
definition of the exponents \eqref{pn},
and the constants $h(T)$ in \eqref{ft} together with
\eqref{jmform}.
\epf

An important particular case of Theorem~\ref{thm:newfus}
and Corollary~\ref{cor:newfuseq} is the case where $\la$
is a partition of $n$. Here the updown tableaux $T$
of shape $\la$ can be regarded as the {\it standard tableaux\/}
which parameterize basis vectors of the corresponding
representation of $\Bc_n(\om)$. It was pointed out in
\cite[Corollary~3.7]{im:fp} that in this
situation all exponents $p_i$ are equal to zero
and the constant $f(T)$ arising in the evaluations \eqref{reev}
depends only on $\la$ and does
not depend on the standard $\la$-tableau $T$.
The new version of the fusion procedure associated with
the rational function \eqref{omubra} possesses the same
property as the next corollary shows. To formulate the result
introduce the {\it content polynomial\/}
\ben
C_{\la}(z)=\prod_{\al\in\la}\big(z+\si(\al)\big),
\een
where $\si(\al)=j-i$ if the box $\al$ of the diagram $\la$
is in row $i$ and column $j$.

\bco\label{cor:symbra}
If $T=(\La_1,\dots,\La_n)$
is an updown $\la$-tableau and $\la=\La_n$ is a partition of $n$,
then the consecutive evaluations give
\ben
\Omega(u_1,\dots,u_n)\big|_{u_1=c_1}\big|_{u_2=c_2}\dots
\big|_{u_n=c_n}=
\frac{2^n\ts C_{\la}(\om/4)\ts H(\la)}{C_{\la}(\om/2)}\ts E_T,
\een
where $H(\la)$ is the product of
the hooks of $\la$.
\eco

\bpf
Arguing as in the proof of Corollary~\ref{cor:newfuseq}
and using the respective calculation of \cite{m:fp}
we observe that the recurrence relation for the constants
$h(T)$ now takes the form
\ben
h(T)=h(U)\ts \frac{H(\la)}{H(\mu)}\ts\frac{\om+4\si_n}{\om+2\si_n},
\een
where $\mu$ is the shape of the tableau $U$ and $\si_n=\si(\al)$
for the box $\al$ occupied by $n$. An obvious induction
completes the proof.
\epf

As the group algebra $\CC[\Sym_n]$ can be regarded as the quotient
of the Brauer algebra $\Bc_n(\om)$ by the ideal generated
by the elements $\ep_1,\dots,\ep_{n-1}$, Corollary~\ref{cor:symbra}
yields a new version of the fusion procedure for the symmetric
group. Moreover, this procedure involves an arbitrary parameter
$\om$ inherited from $\Bc_n(\om)$. To formulate the
corresponding version introduce the rational function
in variables $v_1,\dots,v_n$ with values in $\CC[\Sym_n]$ by
\beql{omvsym}
\Omega_{\om}(v_1,\dots,v_n)=
\prod_{1\leqslant i<j\leqslant n}
\Big(1+\frac{s_{ij}}{v_i+v_j+\om/2}\Big)
\prod_{1\leqslant i<j\leqslant n}
\Big(1-\frac{s_{ij}}{v_i-v_j}\Big)
\eeq
with both products taken in the
lexicographic order on the pairs $(i,j)$.

\bco\label{cor:sym}
Suppose that $\la$ is a partition of $n$ and
$T$ is a standard $\la$-tableau. Let $\si_k=j-i$ if
$k$ occupies the box $(i,j)$ in $T$.
Then the consecutive evaluations give
\ben
\Omega_{\om}(v_1,\dots,v_n)\big|_{v_1=\si_1}\big|_{v_2=\si_2}\dots
\big|_{v_n=\si_n}=
\frac{2^n\ts C_{\la}(\om/4)\ts H(\la)}{C_{\la}(\om/2)}\ts E_T,
\een
where $H(\la)$ is the product of
the hooks of $\la$ and $E_T$ is the primitive idempotent
in $\CC[\Sym_n]$ associated with $T$.
\eco

\bpf
This is immediate from Corollary~\ref{cor:symbra}, where
the relation between the variables is given by
$u_i=v_i+c_1$ for all $i$ and we used the relation
$2\tss c_1-\vk=\om/2$.
\epf

This version of the fusion procedure for the symmetric group
appears to be new. We will discuss its meaning in the context
of evaluation homomorphisms below in Sec.~\ref{sec:tr}.
By taking the limit $\om\to\infty$ we recover
the standard fusion procedure originated in Jucys~\cite{j:yo}
in the form found in \cite{m:fp}; cf. \cite{c:sb} and \cite{n:yc}.

\section{Evaluation homomorphisms}\label{sec:eh}
\setcounter{equation}{0}

Let $G=[g_{ij}]$ be a nonsingular symmetric or skew-symmetric
$N\times N$ matrix with entries in $\CC$ so that
$G^{\tss t}={\pm}\ts G$, where $t$ denotes the standard
matrix transposition. The skew-symmetric case
may occur only if $N$ is even.
Consider the canonical basis $e_1,\dots,e_N$
of the vector space $\CC^N$ and equip it with
the bilinear form
$\langle\ts,\rangle$ by
\beql{formg}
\langle e_i,e_j\rangle=g_{ij}.
\eeq
The classical orthogonal and symplectic groups
$\Or_N$ and $\Spr_N$ consist of the
$N\times N$ matrices
$\hb$ preserving the symmetric
or skew-symmetric form, respectively,
\beql{defosp}
\langle \hb\tss v,\hb\tss w\rangle=\langle v, w\rangle
\qquad\text{for all}\quad v,w\in\CC^N.
\eeq
We will be using
the following convention.
Whenever the double sign $\pm{}$ or $\mp{}$ occurs,
the upper sign will correspond to the symmetric case
and the lower sign to
the skew-symmetric case.
For any $N\times N$ matrix $A$ set
\beql{transG}
A'=G A^t\ts G^{-1}.
\eeq
We denote by $E_{ij}$, $1\leqslant i,j\leqslant N$,
the standard basis vectors of the Lie
algebra $\gl_N$. They satisfy the commutation relations
\beql{comrelgln}
[E_{ij},E_{kl}]=\de_{kj}E_{il}-\de_{il}E_{kj}.
\eeq
For the entries of the inverse matrix of $G$ we will
write $G^{-1}=[\overline{g}_{ij}]$.
Introduce the elements $F_{ij}$
of the Lie algebra $\gl_N$ by the formulas
\beql{fijelem}
F_{ij}=E_{ij}-\sum_{k,\tss l=1}^N
g^{}_{ik}\tss \overline{g}^{}_{lj}\tss E_{lk}.
\eeq
The Lie subalgebra
of $\gl_N$ spanned by the elements $F_{ij}$ is isomorphic to
the orthogonal Lie algebra $\oa_N$ associated with $\Or_N$
in the symmetric case
and to the symplectic Lie algebra $\spa_N$
associated with $\Spr_N$
in the skew-symmetric case.
This Lie algebra will be denoted by $\g_N$.
Common choices of the matrix $G$ in the symmetric case
include the identity matrix $G=1$ and the
antidiagonal matrix $G=[\de_{i,N-j+1}]$. In the skew-symmetric case
with $N=2n$ the entries of $G$ are often chosen in the form
$g_{ij}=\de_{i,2n-j+1}$ for $i=1,\dots,n$ and
$g_{ij}={}-\de_{i,2n-j+1}$ for $i=n+1,\dots,2n$.
In all these cases, the summation in the formula \eqref{fijelem}
reduces to one term.

Consider the endomorphism algebra $\End\CC^N$ and let
$e_{ij}\in\End\CC^N$
be the standard matrix units. We denote by $F$ the $N\times N$ matrix
whose $ij$-th entry is $F_{ij}$. We shall also regard $F$ as the
element
\ben
F=\sum_{i,j=1}^N e_{ij}\ot F_{ij}\in \End\CC^N\ot \U(\g_N).
\een
The definition \eqref{fijelem} can be written
in the matrix form as
\beql{feeg}
F=E-E^{\tss\prime}=E-G\tss E^{\tss t}\tss G^{-1},
\eeq
where
\ben
E=\sum_{i,j=1}^N e_{ij}\ot E_{ij}\in \End\CC^N\ot \U(\gl_N).
\een

Consider the permutation operator
\beql{pmatrix}
P=\sum_{i,j=1}^N e_{ij}\ot e_{ji}\in \End \CC^N\ot\End \CC^N
\eeq
and its partial transpose $P^{\tss t}$ with respect to the first
(or, equivalently, the second) copy of $\End \CC^N$:
\beql{qmatrix}
P^{\tss t}=\sum_{i,j=1}^N e_{ij}\ot e_{ij}.
\eeq
Furthermore, set
\beql{rprimeq}
Q=G_1\tss P^{\tss t}\ts G^{-1}_1=G_2\tss P^{\tss t}\ts G^{-1}_2,
\eeq
where the second equality
follows from the relations $G_1\ts P^{\tss t}=G_2^t\ts P^{\tss t}$
and $P^{\tss t}\ts G_1 =P^{\tss t}\ts G_2^t$.

Note that the operators $P$ and $Q$ satisfy
the relations
\beql{pq}
P^2=1,\qquad PQ=QP={}\pm Q,\qquad Q^2=N\ts Q.
\eeq
The defining relations of the algebra $\U(\gl_N)$ can be written
in a matrix form as
\ben
E_1\ts E_2-E_2\ts E_1=E_1\ts P-P\ts E_1,
\een
where both sides are regarded as elements of the algebra
$\End\CC^N\ot\End\CC^N\ot \U(\gl_N)$ and
\beql{eonetwo}
E_1=\sum_{i,j=1}^N e_{ij}\ot 1\ot E_{ij},\qquad
E_2=\sum_{i,j=1}^N 1\ot e_{ij}\ot E_{ij}.
\eeq
The defining relations of the algebra $\U(\g_N)$
then take the form
\beql{deff}
F_1\ts F_2-F_2\ts F_1=F_1\ts (P-Q)-(P-Q)\ts F_1
\eeq
together with the relation $F+F^{\tss\prime}=0$. The latter
implies
\beql{qf}
Q\ts (F_1+F_2)=(F_1+F_2)\ts Q=0.
\eeq

Set
\beql{kappa}
\kappa=N/2\mp 1.
\eeq
The $R$-matrix $R(u)$ is a rational function in a complex parameter $u$
with values in the tensor product algebra
$\End\CC^N\ot\End\CC^N$ defined by
\beql{Ru}
R(u)=1-\frac{P}{u}+\frac{Q}{u-\kappa}.
\eeq
It is well known by \cite{zz:rf}
that $R(u)$ satisfies the Yang--Baxter equation
\beql{yberep}
R_{12}(u)\ts R_{13}(u+v)\ts R_{23}(v)
=R_{23}(v)\ts R_{13}(u+v)\ts R_{12}(u).
\eeq
Here both sides take values in
$\End\CC^N\ot\End\CC^N\ot\End\CC^N$ and the subscripts indicate
the copies of $\End\CC^N$ so that $R_{12}(u)=R(u)\ot 1$ etc.
Clearly, this is a recast of the properties
of the functions $\rho_{ij}(u)$ defined in \eqref{rij};
we will discuss this relationship
in more detail in Sec.~\ref{sec:tr}.

Following the general approach of \cite{d:qg} and \cite{rtf:ql},
define the {\it extended Yangian\/}
$\X(\g_N)$
as an associative algebra with generators
$t_{ij}^{(r)}$, where $1\leqslant i,j\leqslant N$ and $r=1,2,\dots$,
satisfying certain quadratic relations. In order to write them down,
introduce the formal series
\beql{tiju}
t_{ij}(u)=\de_{ij}+\sum_{r=1}^{\infty}t_{ij}^{(r)}\ts u^{-r}
\in\X(\g_N)[[u^{-1}]]
\eeq
and set
\beql{Tu}
T(u)=\sum_{i,j=1}^N e_{ij}\ot t_{ij}(u)
\in \End\CC^N\ot \X(\g_N)[[u^{-1}]].
\eeq
Consider the algebra
$\End\CC^N\ot\End\CC^N\ot \X(\g_N)[[u^{-1}]]$
and introduce its elements $T_1(u)$ and $T_2(u)$ by
\beql{T1T2}
T_1(u)=\sum_{i,j=1}^N e_{ij}\ot 1\ot t_{ij}(u),\qquad
T_2(u)=\sum_{i,j=1}^N 1\ot e_{ij}\ot t_{ij}(u).
\eeq
The defining relations for the algebra $\X(\g_N)$ can then
be written in the form
\beql{RTT}
R(u-v)\ts T_1(u)\ts T_2(v)=T_2(v)\ts T_1(u)\ts R(u-v).
\eeq
The {\it Yangian\/} $\Y(\g_N)$ is then defined as the quotient
of the extended Yangian $\X(\g_N)$
by the relation
\beql{ttra}
T^{\tss\prime}(u+\ka)\ts T(u)=1.
\eeq
Although the matrix $G$ is implicit in the definitions of
$\X(\g_N)$ and $\Y(\g_N)$, the respective
algebras associated with two different
nonsingular symmetric (or skew-symmetric) $N\times N$ matrices
$G$ and $\wt G$ are isomorphic to each other so that
each of them depends only on the Lie algebra $\g_N$.
Drinfeld defined the Yangians associated with
simple Lie algebras by using
different presentations; see \cite{d:ha, d:qg}.
Explicit form of the defining relations for some standard
choices of antidiagonal symmetric and skew-symmetric matrices $G$
in terms of the series $t_{ij}(u)$
can be found in \cite{aacfr:rp}, where the equivalence
of the $R$-matrix presentation and the Drinfeld presentation
is explained; see also \cite{amr:rp}.
The Yangian $\Y(\g_N)$ carries a Hopf algebra structure
defined in a standard way; see e.g. {\it loc. cit.}
for explicit formulas.

Following \cite{s:bc}, use a reflection type equation
to define algebras associated to $\g_N$;
see also \cite{aacfr:rp}, \cite{mr:rr}.

\bde\label{def:refeq}
The {\it reflection algebra\/} $\Br(\g_N)$
corresponding to the Lie algebra $\g_N$ is defined as the
associative algebra generated by elements $s_{ij}^{(r)}$
with $1\leqslant i,j\leqslant N$ and $r\geqslant 1$
subject to defining relations written in terms of the
generating series
\ben
s_{ij}(u)=\de_{ij}+\sum_{r=1}^{\infty}s_{ij}^{(r)}\ts u^{-r}
\een
as follows. Introduce the matrix
\beql{Su}
S(u)=\sum_{i,j=1}^N e_{ij}\ot s_{ij}(u)
\in \End\CC^N\ot \Br(\g_N)[[u^{-1}]].
\eeq
Then the defining relations have the form
of the {\it reflection equation\/}
\beql{refeq}
R(u-v)\ts S_1(u)\ts R(u+v)\ts S_2(v)=S_2(v)\ts
R(u+v)\ts S_1(u)\ts R(u-v)
\eeq
and the {\it unitary condition\/}
\beql{unita}
S(u)\ts S(-u)=1,
\eeq
where we use the matrix notation as in \eqref{RTT}.
\qed
\ede

\bpr\label{prop:subalg}
The mapping
\beql{hom}
S(u)\mapsto T(-u)^{-1}\ts T(u)
\eeq
defines an algebra homomorphism $\Br(\g_N)\to\X(\g_N)$.
\epr

\bpf
It is obvious that \eqref{unita} holds when $S(u)$ is replaced
by the matrix $T(-u)^{-1}\ts T(u)$.
Checking that the image satisfies
\eqref{refeq} is standard and relies on \eqref{RTT}
and the relation
\ben
T_1(u)\ts R(u+v)\ts T_2(-v)^{-1}=T_2(-v)^{-1}\ts R(u+v)\ts T_1(u)
\een
implied by \eqref{RTT}.
\epf

The Yangian
$\Y(\g_N)$ contains the universal enveloping algebra
$\U(\g_N)$ as a subalgebra. However, as shown by
Drinfeld~\cite{d:ha}, in contrast with the
Yangian for the Lie algebra $\gl_N$,
there is no homomorphism $\Y(\g_N)\to \U(\g_N)$, identical
on the subalgebra $\U(\g_N)$. Nevertheless, as the following
Theorem~\ref{thm:eval} and Corollary~\ref{cor:embed} show,
such a homomorphism from the reflection algebra does exist.
This result is suggested by the solutions of the reflection
equation associated with
the affine Birman--Murakami--Wenzl algebras
found in \cite{io:bs}.

\bth\label{thm:eval}
The mapping
\beql{evalya}
S(u)\mapsto \frac{u+F-N/4}{u-F+N/4}
=1+2\ts\sum_{r=1}^{\infty}
\Big(F-\frac{N}{4}\Big)^{\tss r}\ts u^{-r}
\eeq
defines an algebra homomorphism $\Br(\g_N)\to\U(\g_N)$.
\eth

\bpf
Set $\Fc(u)=u+F$. To make our formulas more readable, we will
be using the notation $a=N/4$ throughout the proof.
We will demonstrate that
both relations \eqref{refeq} and
\eqref{unita} will hold when the matrix $S(u)$ is replaced by
the product $\Fc(u-a)\Fc(-u-a)^{-1}$. This is obvious for \eqref{unita},
so we come to proving the identity
\begin{multline}\label{refeqfu}
R(u-v)\ts \Fc_1(u-a)\ts \Fc_1(-u-a)^{-1}
\ts R(u+v)\ts \Fc_2(v-a)\ts \Fc_2(-v-a)^{-1}\\
{}=\Fc_2(v-a)\ts \Fc_2(-v-a)^{-1}\ts
R(u+v)\ts \Fc_1(u-a)\ts \Fc_1(-u-a)^{-1}\ts R(u-v).
\end{multline}

\ble\label{lem:rff}
We have the identity
\beql{rff}
R(u-v)\ts \Fc_1(u)\ts \Fc_2(v)
-\Fc_2(v)\ts \Fc_1(u)\ts R(u-v)=\frac{1}{u-v-\ka}\ts U,
\eeq
where
\ben
U=Q\ts (F_1+\ka)\ts F_2-F_2\ts (F_1+\ka)\ts Q.
\een
\ele

\bpf
The left hand side of \eqref{rff} simplifies to
\ben
F_1\ts F_2-F_2\ts F_1+P\ts F_1-F_1\ts P
+\frac{1}{u-v-\ka}\Big(Q\ts F_1\ts F_2-F_2\ts F_1\ts Q
+(u-v)\big(Q\ts F_2-F_2\ts Q\big)\Big),
\een
where we used the relations $P\ts F_1=F_2\ts P$ and \eqref{qf}.
Now use \eqref{deff} to write
\ben
F_1\ts F_2-F_2\ts F_1+P\ts F_1-F_1\ts P
=Q\ts F_1-F_1\ts Q={}-Q\ts F_2+F_2\ts Q
\een
which gives the required relation.
\epf

We will now establish some simple properties of the element $U$
to be used below.

\ble\label{lem:utwoone}
We have the relations
\begin{align}
\label{puup}
P\ts U=U P&={}\pm{} U,\\
\label{quuq}
Q\ts U+U\ts Q&=N\ts U,\\
\label{ffuq}
(F_1+F_2)\ts U&=U\ts Q.
\end{align}
\ele

\bpf
We have
\ben
P\ts U=P\ts Q\ts (F_1+\ka)\ts F_2-P\ts F_2\ts (F_1+\ka)\ts Q
={}\pm Q\ts (F_1+\ka)\ts F_2\mp\ts F_1\ts (F_2+\ka)\ts Q.
\een
Applying \eqref{pq} and \eqref{deff}, we get
\ben
\bal
F_1\ts F_2\ts Q&=\big(F_2\ts F_1+F_1\ts(P-Q)-(P-Q)\ts F_1\big)\ts Q\\
{}&= F_2\ts F_1\ts Q-(N\mp 1)\ts F_1\ts Q\mp F_2\ts Q
=F_2\ts F_1\ts Q-2\ts\ka\ts F_1\ts Q,
\eal
\een
where we also used \eqref{qf} and
the easily verified relation $Q\ts F_1\ts Q=0$.
This implies $P\ts U=\pm U$. The second relation in \eqref{puup}
is verified by essentially the same calculation.

Relation \eqref{quuq} is immediate from the identity
$Q\tss F_1\tss F_2\tss Q=Q\tss F_2\tss F_1\tss Q$. To prove
\eqref{ffuq} use \eqref{qf} to write
\ben
\bal
(F_1+F_2)\ts U&=-(F_1+F_2)\ts F_2\ts (F_1+\ka)\ts Q
=-[F_1+F_2,F_2\ts(F_1+\ka)]\ts Q\\
{}&=-[F_1,F_2]\ts (F_1+\ka)\ts Q
+F_2\ts [F_1,F_2]\ts Q.
\eal
\een
Due to \eqref{deff} we may replace $[F_1,F_2]$ by $[F_1,P-Q]$
and complete the calculation as in the proof of \eqref{puup}
to show that the expression coincides with $U\ts Q$.
\epf

Multiplying both sides of \eqref{rff} by $\Fc_2(v)^{-1}$ from the left
and the right, we come to
the relation
\beql{rffi}
\Fc_1(u)\ts R(u-v)\ts \Fc_2(v)^{-1}
-\Fc_2(v)^{-1}\ts R(u-v)\ts \Fc_1(u)=
-\frac{1}{u-v-\ka}\ts \Fc_2(v)^{-1}\ts U\ts \Fc_2(v)^{-1}.
\eeq
By conjugating its both sides by $P$ and
using Lemma~\ref{lem:utwoone}, we get
\beql{rffin}
\Fc_1(v)^{-1}\ts R(u-v)\ts \Fc_2(u)
-\Fc_2(u)\ts R(u-v)\ts \Fc_1(v)^{-1}=
\frac{1}{u-v-\ka}\ts \Fc_1(v)^{-1}\ts U\ts \Fc_1(v)^{-1}.
\eeq
We will need one more relation obtained from \eqref{rffin}
by multiplying both sides by $\Fc_2(u)^{-1}$ from the left
and the right:
\begin{multline}\label{rffinv}
R(u-v)\ts \Fc_1(v)^{-1}\ts \Fc_2(u)^{-1}
-\Fc_2(u)^{-1}\ts \Fc_1(v)^{-1}\ts R(u-v)\\
={}
-\frac{1}{u-v-\ka}\ts \Fc_2(u)^{-1}\ts
\Fc_1(v)^{-1}\ts U\ts \Fc_1(v)^{-1}
\ts \Fc_2(u)^{-1}.
\end{multline}

Our strategy in proving \eqref{refeqfu} is to use
Lemma~\ref{lem:rff} and relations \eqref{rffi}--\eqref{rffinv}
to transform the expression on the left hand side
to that on the right hand side and then to show that the
additional terms arising in the process will add up to zero.
We start by applying \eqref{rffin} with
the substitution $u\mapsto v-a$ and $v\mapsto -u-a$.
This yields the expression
\ben
R(u-v)\ts \Fc_1(u-a)\ts \Fc_2(v-a)
\ts R(u+v)\ts \Fc_1(-u-a)^{-1}\ts \Fc_2(-v-a)^{-1}
\een
together with the additional term
\beql{termone}
\frac{1}{u+v-\ka}\ts R(u-v)\ts \Fc_1(u-a)\ts
\Fc_1(-u-a)^{-1}\ts U\ts \Fc_1(-u-a)^{-1}\ts \Fc_2(-v-a)^{-1}.
\eeq
Now use \eqref{rff} with the substitution
$u\mapsto u-a$ and $v\mapsto v-a$ to get
\ben
\Fc_2(v-a)\ts \Fc_1(u-a)\ts R(u-v)
\ts R(u+v)\ts \Fc_1(-u-a)^{-1}\ts \Fc_2(-v-a)^{-1}
\een
with the additional term
\beql{termtwo}
\frac{1}{u-v-\ka}\ts U\ts R(u+v)\ts \Fc_1(-u-a)^{-1}\ts \Fc_2(-v-a)^{-1}.
\eeq
Furthermore, since
$R(u-v)\ts R(u+v)=R(u+v)\ts R(u-v)$, applying
\eqref{rffinv} with the substitution
$u\mapsto -v-a$ and $v\mapsto -u-a$ we get the expression
\ben
\Fc_2(v-a)\ts \Fc_1(u-a)
\ts R(u+v)\ts \Fc_2(-v-a)^{-1}\ts \Fc_1(-u-a)^{-1}\ts R(u-v)
\een
together with the additional term
\begin{multline}\label{termthree}
{}-\frac{1}{u-v-\ka}\ts
\Fc_2(v-a)\ts \Fc_1(u-a)\ts R(u+v)\ts \Fc_2(-v-a)^{-1}\\
{}\times\ts \Fc_1(-u-a)^{-1}\ts U
\ts \Fc_1(-u-a)^{-1}\ts \Fc_2(-v-a)^{-1}.
\end{multline}
Finally, the application of \eqref{rffi} with the substitution
$u\mapsto u-a$ and $v\mapsto -v-a$ yields the right hand
side of \eqref{refeqfu} with the additional term
\beql{termfour}
{}-\frac{1}{u+v-\ka}\ts \Fc_2(v-a)\ts \Fc_2(-v-a)^{-1}\ts U
\ts \Fc_2(-v-a)^{-1}\ts \Fc_1(-u-a)^{-1}\ts R(u-v).
\eeq

We need to show that the sum of the four additional terms
\eqref{termone}--\eqref{termfour} is zero. To do this, note
that the sum of the terms \eqref{termthree} and \eqref{termfour}
equals the sum of the terms
\begin{multline}\label{termthreenew}
{}-\frac{1}{u-v-\ka}\ts
\Fc_2(v-a)\ts \Fc_2(-v-a)^{-1}\ts  R(u+v)\ts \Fc_1(u-a)\\
{}\times\ts \Fc_1(-u-a)^{-1}\ts U
\ts \Fc_1(-u-a)^{-1}\ts \Fc_2(-v-a)^{-1}
\end{multline}
and
\beql{termfournew}
{}-\frac{1}{u+v-\ka}\ts \Fc_2(v-a)\ts \Fc_2(-v-a)^{-1}\ts U
\ts R(u-v)\ts \Fc_1(-u-a)^{-1}\ts \Fc_2(-v-a)^{-1}.
\eeq
Indeed, this follows by applying \eqref{rffi} with the substitution
$u\mapsto u-a$ and $v\mapsto -v-a$ to \eqref{termthree},
and applying \eqref{rffinv}
with the substitution
$u\mapsto -v-a$ and $v\mapsto -u-a$ to \eqref{termfour}.
Clearly, the new additional terms arising in these transformations
cancel.

Now multiply each of the expressions in \eqref{termone},
\eqref{termtwo}, \eqref{termthreenew} and \eqref{termfournew}
by the product $\Fc_2(-v-a)\ts \Fc_1(-u-a)$ from the right and by
$(u-v-\ka)(u+v-\ka)\ts \Fc_2(-v-a)$ from the left.
Adding up the resulting expressions we get
\ben
\bal
&(u-v-\ka)\ts \Fc_2(-v-a)\ts R(u-v)\ts \Fc_1(u-a)\ts
\Fc_1(-u-a)^{-1}\ts U\\
{}+{}&(u+v-\ka)\ts \Fc_2(-v-a)\ts U\ts R(u+v)\\
{}-{}&(u+v-\ka)\ts
\Fc_2(v-a)\ts  R(u+v)\ts \Fc_1(u-a)\ts \Fc_1(-u-a)^{-1}\ts U\\
{}-{}&(u-v-\ka)\ts \Fc_2(v-a)\ts U\ts R(u-v).
\eal
\een
Using the definition \eqref{Ru} of $R(u)$ we bring this
sum to the form
\ben
\bal
2\tss v{}&\ts\Big(P-Q-\Fc_2(u-\ka-a)\Big)\Fc_1(u-a)\ts
\Fc_1(-u-a)^{-1}\ts U\\
{}&{}+2\tss v\Big(U(P-Q)+\Fc_2(-u+\ka-a)\ts U\Big)\\
{}&{}+\frac{2v\tss\ka}{u^2-v^2}\ts \Big(\Fc_2(-u-a)
\ts P\ts \Fc_1(u-a)\ts
\Fc_1(-u-a)^{-1}\ts  U - \Fc_2(u-a)\ts U P\Big).
\eal
\een
By Lemma~\ref{rprimeq} we have
\ben
P\ts \Fc_1(u-a)\ts
\Fc_1(-u-a)^{-1}\ts  U=\Fc_2(u-a)\ts
\Fc_2(-u-a)^{-1}\ts  U\ts P
\een
so that the last summand vanishes. Thus, it remains to
verify that
\begin{multline}\label{pqfze}
\Big(P-Q-\Fc_2(u-\ka-a)\Big)\Fc_1(u-a)\ts
\Fc_1(-u-a)^{-1}\ts U\\
{}+{}U(P-Q)+\Fc_2(-u+\ka-a)\ts U=0.
\end{multline}
However, by \eqref{deff}, the term $P-Q-\Fc_2(u-\ka-a)$
commutes with $F_1$. Therefore, when multiplied from the left
by $\Fc_1(-u-a)$, the left hand side of \eqref{pqfze} becomes
\ben
\Fc_1(u-a)\ts\Big(P-Q-\Fc_2(u-\ka-a)\Big)\ts U
+\Fc_1(-u-a)\ts \Big(U(P-Q)+\Fc_2(-u+\ka-a)\ts U\Big),
\een
which simplifies to
\ben
u\Big(N\ts U+(P-Q)\ts U-U\ts(P-Q)-2(F_1+F_2)\ts U\Big)
+(F_1-a)\Big((P-Q)\ts U+U\ts(P-Q)-2\ka\ts U\Big).
\een
Lemma~\ref{lem:utwoone} implies that this expression
is zero thus proving \eqref{pqfze} and the theorem.
\epf

\bre\label{rem:affbra}
The evaluation homomorphism of Theorem~\ref{thm:eval}
can be ``lifted" to
the {\it affine Brauer algebra\/} defined in \cite{n:yo} (under
the name {\it degenerate affine Wenzl algebra\/}).
The precise statement will be given in \cite{imo:fp} in the
context of the affine Birman--Murakami--Wenzl algebras
in the spirit of the approach to the representation theory
of these algebras developed in \cite{io:jm}.
\qed
\ere

In the following we denote by $s^{\tss\prime}_{ij}(u)$ the entries
of the matrix $S^{\tss\prime}(u)=G\tss S^{\tss t}(u)\ts G^{-1}$,
and write
\ben
s^{\tss\prime}_{ij}(u)=\de_{ij}
+\sum_{r=1}^{\infty}s^{\tss\prime\tss(r)}_{ij}\ts u^{-r}.
\een

\bco\label{cor:embed}
The assignment
\beql{embedlie}
F_{ij}\mapsto \frac14\ts\Big(s_{ij}^{(1)}
-s^{\tss\prime\tss(1)}_{ij}\Big)
\eeq
defines an embedding $\U(\g_N)\hra\Br(\g_N)$. Moreover,
the homomorphism \eqref{evalya} is identical
on $\U(\g_N)$.
\eco

\bpf
Write \eqref{embedlie} in a matrix form
\beql{matemb}
F\mapsto \frac14\ts\big(S-S^{\tss\prime}\big),
\eeq
with $S=[s_{ij}^{(1)}]$ and
$S^{\tss\prime}=[s^{\tss\prime\tss(1)}_{ij}]$, and verify that
this assignment defines a homomorphism. Obviously,
the relation $F+F'=0$ is preserved by the assignment so we
are left to show that \eqref{deff} is preserved as well.
Expand the rational functions in $u$ and $v$ involved in
\eqref{refeq} into series in $v^{-1}$ by
\ben
\frac{1}{v-a}=v^{-1}+a\tss v^{-2}+\dots
\een
and compare the coefficients of $u^{-1}v^{-1}$ on both sides.
This yields
\ben
S_1\ts S_2-S_2\ts S_1=2\ts S_1\ts(P-Q)-2\ts (P-Q)\ts S_1.
\een
Now apply partial transpositions with respect to the first
and the second copies of $\End\CC^N$ to get
\ben
\bal
S^{\tss\prime}_1\ts S_2-S_2\ts S^{\tss\prime}_1&
=2\ts S^{\tss\prime}_1\ts(P-Q)-2\ts (P-Q)\ts S^{\tss\prime}_1,\\
S_1\ts S^{\tss\prime}_2-S^{\tss\prime}_2\ts S_1&
=2\ts (P-Q)\ts S_1-2\ts S_1\ts(P-Q),\\
S^{\tss\prime}_1\ts S^{\tss\prime}_2
-S^{\tss\prime}_2\ts S^{\tss\prime}_1&
=2\ts (P-Q)\ts S^{\tss\prime}_1-2\ts S^{\tss\prime}_1\ts(P-Q).
\eal
\een
This implies that \eqref{deff} holds under
the assignment \eqref{matemb} and thus the latter
defines a homomorphism.

Furthermore, under the homomorphism \eqref{evalya}
we have $S\mapsto 2\ts F-N/2$ and so its composition
with \eqref{matemb} (applying \eqref{matemb} first)
is the identity map on $\U(\g_N)$.
Hence the kernel of the homomorphism \eqref{matemb} is zero.
\epf

\bre\label{rem:hom}
Due to relation \eqref{ttra}, the homomorphism \eqref{hom}
provides a homomorphism $\Br(\g_N)\to\Y(\g_N)$ which
can be written as
\beql{twhom}
S(u)\to T^{\tss\prime}(-u+\ka)\ts T(u).
\eeq
As pointed out in \cite{aacfr:rp},
this brings up a connection of the reflection algebra $\Br(\g_N)$
with the Olshanski {\it twisted Yangians\/} \cite{o:ty},
see also \cite[Ch.~2]{m:yc}. The homomorphism
\eqref{twhom} is not injective and its kernel
can be described by a symmetry-type relations; cf. {\it loc. cit.}
For instance,
in the case $G=1$ the elements
$s_{ii}^{(1)}$ belong to the kernel. However, $s_{ii}^{(1)}\ne 0$
in $\Br(\g_N)$, as follows from Theorem~\ref{thm:eval}.
In other words, if we define the twisted
Yangian $\Y^{\rm tw}(\g_N)$
as the quotient of $\Br(\g_N)$ by the kernel
of the homomorphism \eqref{twhom},
then the evaluation homomorphism of Theorem~\ref{thm:eval}
will not factor through the natural epimorphism
$\Br(\g_N)\to \Y^{\rm tw}(\g_N)$.
\qed
\ere

\bre\label{rem:super}
Most of the results of this paper
admit natural super-analogues where the Lie algebra $\g_N$
is replaced by the orthosymplectic Lie superalgebra $\osp_{m|2n}$.
Indeed, our calculations are performed in
a matrix language so that to apply the same approach in the super
case we only need to set up appropriate matrix notation
taking care of the sign rules. Here we will indicate
the necessary changes to be made in the notation
and formulas. The $\ZZ_2$-degree (or parity) of
the basis elements $E_{ij}$ of the
Lie superalgebra $\gl_{m|2n}$ is given by $\deg(E_{ij})=\bi+\bj$,
where $\bi=0$ for $1\leqslant i\leqslant m$ and
$\bi=1$ for $m+1\leqslant i\leqslant m+2n$.
The commutation relations
in $\gl_{m|2n}$ have the form
\ben
[E_{ij}, E_{kl}]=\de_{kj}E_{il}-\de_{il}E_{kj}
(-1)^{(\bi+\bj)(\bk+\bl)},
\een
where the square brackets denote the super-commutator.
We will work with square matrices $A=[A_{ij}]$
of size $m+2n$. Any such matrix with entries in
a superalgebra $\Ac$ will be identified with
the element
\ben
\sum_{i,j=1}^{m+2n}e_{ij}\ot A_{ij}(-1)^{\bi\bj+\bj}
\in \End\CC^{m|2n}\ot\Ac.
\een
The signs here are necessary to preserve matrix
multiplication rules as we work with graded tensor products of
superalgebras.
Fix a block-diagonal matrix $G$
whose upper-left $m\times m$ block is a nonsingular
symmetric matrix while the lower-right $2n\times 2n$
block is nonsingular skew-symmetric matrix.
Given a matrix $A$ we define its transpose associated with $G$
by the formula \eqref{transG}, where $t$ now denotes
a standard matrix super-transposition,
$(A^t)_{ij}=A_{ji}(-1)^{\bi\bj+\bi}$. Note that in contrast
with the standard super-transposition,
the transposition defined by \eqref{transG} is involutive.
We will also regard $t$ as a linear map
\beql{suptra}
t:\End\CC^{m|2n}\to \End\CC^{m|2n}, \qquad
e_{ij}\mapsto e_{ji}(-1)^{\bi\bj+\bj}.
\eeq
In the case of multiple tensor products of the superalgebras
$\End\CC^{m|2n}$ we will indicate by $t_a$ the map \eqref{suptra}
acting on the $a$-th copy.
Introduce the square matrix $E=[E_{ij}(-1)^{\bj}]$
of size $m+2n$ and define the elements $F_{ij}$ of $\gl_{m|2n}$
by the formula \eqref{feeg}, where $F$ is the matrix
$F=[F_{ij}(-1)^{\bj}]$. The Lie superalgebra
spanned by the elements $F_{ij}$ is isomorphic to
the orthosymplectic Lie superalgebra $\osp_{m|2n}$.
The matrix form of the commutation relations in $\osp_{m|2n}$
coincides with \eqref{deff} together with $F+F^{\tss\prime}=0$;
the latter implies \eqref{qf}. This time
both sides of \eqref{deff} are elements of the superalgebra
\ben
\End\CC^{m|2n}\ot \End\CC^{m|2n}\ot \U(\osp_{m|2n}),
\een
while the definitions of the operators $P$ and $Q$ are
modified by
\ben
P=\sum_{i,j=1}^{m+2n} e_{ij}\ot e_{ji}(-1)^{\bj}
\in \End\CC^{m|2n}\ot \End\CC^{m|2n}
\een
and
\ben
Q=G_1\tss P^{\tss t_1}\ts G_1^{-1}=G_2\tss P^{\tss t_2}\ts G_2^{-1}.
\een
Instead of \eqref{pq} we have
\ben
P^2=1,\qquad PQ=QP=Q,\qquad Q^2=(m-2n)\ts Q.
\een
The $R$-matrix $R(u)$ has the form \eqref{Ru} with
$\ka=m/2-n-1$ which leads to the definitions of the
Yangian and the reflection algebra $\Br(\osp_{m|2n})$
associated with $\osp_{m|2n}$ as in Definition~\ref{def:refeq};
see \cite{aacfr:rp}. The evaluation homomorphism
$\Br(\osp_{m|2n})\to\U(\osp_{m|2n})$ provided by
Theorem~\ref{thm:eval} takes the same form \eqref{evalya},
where $N$ should be replaced by $m-2n$:
\ben
S(u)\mapsto \frac{u+F-(m-2n)/4}{u-F+(m-2n)/4}.
\een
The embedding $\U(\osp_{m|2n})\hra\Br(\osp_{m|2n})$ is defined
by the same formula \eqref{matemb}.
Note also that the
parameter $\om$ of the Brauer algebra $\Bc_n(\om)$
should be evaluated at $\om=m-2n$ in the context
of the super-version of the centralizer results
involving $\osp_{m|2n}$;
see Sec.~\ref{sec:tr}.
\qed
\ere

\section{Tensor representations of
reflection algebras}\label{sec:tr}
\setcounter{equation}{0}

The action of the Lie algebra $\gl_N$
in the vector representation $\CC^N$ is given by the rule
$E_{ij}\mapsto e_{ij}$ which corresponds to the natural action
of the general linear group ${\rm GL}_N(\CC)$ on $\CC^N$
by left multiplication.
By restricting this action to the subgroups
of orthogonal and symplectic matrices
$\Or_N$ and $\Spr_N$ ($N$ is even for the latter)
we make $\CC^N$ and the
tensor product space
\beql{cntenpr}
\CC^N\ot\CC^N\ot\dots\ot\CC^N,\qquad n\ \ \text{factors},
\eeq
into representations of $\Or_N$ and $\Spr_N$.
To write the images of elements of the groups,
introduce an extra copy of the vector space $\CC^N$ and consider
the endomorphism algebra
\beql{endspp}
\End\CC^N\ot\End(\CC^N)^{\ot n},
\eeq
labeling the tensor factors with the numbers
$0,1,\dots,n$. Writing a group element $\hb$ as
$
\hb=\sum_{i,j=1}^N e_{ij}\ot \hb_{ij},
$
for the image of $\hb$ we have
$
\hb\mapsto \hb_1\dots\hb_n
$
with the notation similar to \eqref{eonetwo} and \eqref{T1T2}.
Due to the work of
Brauer~\cite{b:aw}, the centralizer of this action
in the endomorphism algebra of \eqref{cntenpr} is the homomorphic
image of the algebra $\Bc_n(\om)$
with an appropriately specified
value of the parameter $\om$. Namely, $\om=N$
in the orthogonal case, and
the action of the Brauer algebra $\Bc_n(N)$
in the vector space \eqref{cntenpr}
is given by
\beql{braact}
s_{ij}\mapsto P_{ij},\qquad \ep_{ij}\mapsto Q_{ij},\qquad
i< j.
\eeq
In the symplectic case, $\om=-N$
and the action of the Brauer algebra $\Bc_n(-N)$
is given by
\beql{braactsy}
s_{ij}\mapsto -\ts P_{ij},\qquad \ep_{ij}\mapsto -\ts Q_{ij},\qquad
i< j;
\eeq
see \cite{b:aw}, \cite{w:sb},
where the operators $P_{ij}$ and $Q_{ij}$ on
the vector space \eqref{cntenpr} act as
the respective operators $P$ and $Q$ defined in \eqref{pmatrix}
and \eqref{rprimeq} on the tensor product on the $i$-th and $j$-th
copies of $\CC^N$ and act as the identity operators on each
of the remaining copies.
Furthermore, in the orthogonal case
the vector space \eqref{cntenpr} is decomposed
as
\ben
(\CC^N)^{\ot n}\cong \bigoplus_{f=0}^{\lfloor n/2\rfloor}
\bigoplus_{\overset{\scriptstyle
\la\vdash n-2f}{\la'_1+\la'_2\leqslant N}}V_{\la}\ot L(\la),
\een
where $V_{\la}$ and $L(\la)$
are the respective irreducible representations of $\Bc_n(N)$ and $\Or_N$
associated with the diagram $\la$, and we denote
by $\la'$ the conjugate
diagram so that
$\la'_j$ is the number
of boxes in the column $j$ of $\la$; see \cite{w:cg}.
Similarly, in the symplectic case
\ben
(\CC^N)^{\ot n}\cong \bigoplus_{f=0}^{\lfloor n/2\rfloor}
\bigoplus_{\overset{\scriptstyle
\la\vdash n-2f}{2\la'_1\leqslant N}}V_{\la'}\ot L(\la),
\een
where $V_{\la'}$ and $L(\la)$
are the respective irreducible representations of
$\Bc_n(-N)$ and $\Spr_N$
associated with $\la'$ and $\la$; see {\it loc. cit.}

From now on we will impose the condition on the parameters,
\beql{condnn}
N\geqslant 2n+ 1\Fand  N\geqslant 2n
\eeq
in the orthogonal and symplectic case, respectively.
It implies, in particular, that
the Brauer algebras $\Bc_n(N)$ and $\Bc_n(-N)$ are semisimple;
see e.g. \cite{r:cs} where a criterion of semisimplicity of
$\Bc_n(\om)$ is given.
Explicit projections of $(\CC^N)^{\ot n}$
on the irreducible representations of the orthogonal
and symplectic groups are provided by the idempotents
of the corresponding Brauer algebra.
More precisely, suppose that $\Uc=(\La_1,\dots,\La_n)$
is an updown tableau of shape $\la=\La_n$.
Consider the corresponding idempotent $E_{\Uc}$,
the sequence of contents $(c_1,\dots,c_n)$ associated with $\Uc$
and the subspace
\beql{sublu}
L_{\Uc}=E_{\Uc}(\CC^N)^{\ot n}.
\eeq
The condition \eqref{condnn} ensures that $E_{\Uc}$ can be
defined by \eqref{murphyfo} with non-vanishing denominators.
Under the action \eqref{braact} of the Brauer
algebra $\Bc_n(N)$ the subspace
is nonzero if the number of boxes in the first two
columns of $\la$ does not exceed $N$, which follows automatically
from \eqref{condnn}.
In this case
$L_{\Uc}$ is
an irreducible representation of $\Or_N$ isomorphic to $L(\la)$.
Similarly, under the action \eqref{braactsy} of the Brauer
algebra $\Bc_n(-N)$ the subspace \eqref{sublu} is nonzero
if the number of boxes in the first
column of $\la$ does not exceed $N/2$;
this holds automatically
by \eqref{condnn}.
In this case
$L_{\Uc}$ is
an irreducible representation of $\Spr_N$ isomorphic to $L(\la)$.

The vector space $L_{\Uc}$ carries a representation
of the corresponding Lie algebra $\g_N=\oa_N$ or $\g_N=\spa_N$.
Employing the evaluation homomorphism \eqref{evalya},
we can equip $L_{\Uc}$ with the action
of the reflection algebra $\Br(\g_N)$ associated with
the Lie algebra $\g_N$.

On the other hand, the vector space $(\CC^N)^{\ot n}$ carries
a representation of the extended Yangian $\X(\g_N)$ defined
by the assignment
\beql{yarep}
T(u)\mapsto a(u)\ts R_{01}(u-z_1)\ts
\dots R_{0n}(u-z_n),
\eeq
where $z_1,\dots,z_n$ are fixed complex numbers and
$a(u)\in 1+\CC[[u^{-1}]]\ts u^{-1}$ is a fixed
formal series.
Here we regard the image of $T(u)$ as a formal series in $u^{-1}$
with coefficients in the algebra \eqref{endspp}.
Indeed, verifying that \eqref{yarep} defines
a representation of the Yangian amounts to checking that the image
of $T(u)$ satisfies the defining relations \eqref{RTT},
which is straightforward. Furthermore, using
Proposition~\ref{prop:subalg}, we can equip
the vector space $(\CC^N)^{\ot n}$ with the action
of the reflection algebra $\Br(\g_N)$ by
\beql{sua}
S(u)\mapsto a(u)\tss a(-u)^{-1}\ts R_{0n}(-u-z_n)^{-1}\ts
\dots R_{01}(-u-z_1)^{-1}\ts R_{01}(u-z_1)\ts
\dots R_{0n}(u-z_n).
\eeq

We keep using the double sign convention as in Sec.~\ref{sec:eh}
so that the upper sign is taken in the orthogonal case
and the lower sign in the symplectic case.

\bpr\label{prop:invco}
The subspace $L_{\Uc}$ of the vector space \eqref{cntenpr}
is invariant under the action of $\Br(\g_N)$
given by
\begin{multline}
S(u)\mapsto \frac{u-N/4}{u+N/4}\ts
\ts R_{0n}(-u\mp c_n+\ka/2)^{-1}\ts
\dots R_{01}(-u\mp c_1+\ka/2)^{-1}\\
{}\times R_{01}(u\mp c_1+\ka/2)\ts
\dots R_{0n}(u\mp c_n+\ka/2).
\non
\end{multline}
Moreover, the representation of $\Br(\g_N)$ on $L_{\Uc}$
is isomorphic to the evaluation module
$L_{\Uc}$ obtained via the homomorphism \eqref{evalya}.
\epr

\bpf
Consider the anti-automorphism of the Brauer algebra $\Bc_n(\om)$
identical on the generators. Apply this
anti-automorphism to the identity \eqref{jmm},
then replace $u$ by $u+\vk/2$. Together with \eqref{unitarm}
this leads to the following relation in the orthogonal
case:
\begin{multline}
R_{0n}(-u-c_n+\ka/2)^{-1}\ts
\dots R_{01}(-u-c_1+\ka/2)^{-1}\\
{}\times R_{01}(u-c_1+\ka/2)\ts
\dots R_{0n}(u-c_n+\ka/2)
\ts E_{\Uc}
=\frac{u+N/4}{u-N/4}\ts E_{\Uc}\ts \frac{u+X_0-N/4}{u-X_0+N/4},
\non
\end{multline}
where $X_0=Q_{01}+\dots+Q_{0,n}-P_{01}-\dots-P_{0,n}$
so that the parameter $n$ in \eqref{jmm} is replaced by $n+1$
and $m$ is replaced by $0$. We have used the fact that
the image of $\rho_{ij}(u)$ coincides with $R_{ij}(u)$
under the specialization \eqref{braact} with $\om=N$.

To get the symplectic counterpart
of the relation, note that under
the specialization \eqref{braactsy} with $\om=-N$
the image of $\rho_{ij}(u)$ coincides with $R_{ij}(-u)$.
Before applying the above argument, we invert all factors of $E_U$
which occur in \eqref{jmm} so that in the symplectic case
we get
\begin{multline}
R_{0n}(-u+c_n+\ka/2)^{-1}\ts
\dots R_{01}(-u+c_1+\ka/2)^{-1}\\
{}\times R_{01}(u+c_1+\ka/2)\ts
\dots R_{0n}(u+c_n+\ka/2)
\ts E_{\Uc}
=\frac{u+N/4}{u-N/4}\ts E_{\Uc}\ts \frac{u+X_0-N/4}{u-X_0+N/4}.
\non
\end{multline}
This implies the first part of the proposition. The second part
follows from the observation that $X_0$ coincides with
the image of the element $F$ in the representation
\eqref{cntenpr} so that
the two actions of the reflection algebra $\Br(\g_N)$
on $L_{\Uc}$ coincide.
\epf

We conclude with an interpretation of the new fusion procedure
for the symmetric group provided by Corollary~\ref{cor:sym},
from the viewpoint of the representation theory of the Yangians.
Recall that the {\it Yangian for\/} $\gl_N$ is
an associative algebra $\Y(\gl_N)$ with generators
$t_{ij}^{(r)}$, where $1\leqslant i,j\leqslant N$ and $r=1,2,\dots$,
satisfying certain quadratic relations. They are written
with the use of the generating functions \eqref{tiju}
and have the form of the $RTT$ relation \eqref{RTT},
where instead of the $R$-matrix \eqref{Ru} we take the
{\it Yang $R$-matrix\/}
\beql{YRu}
R(u)=1-\frac{P}{u};
\eeq
see e.g. \cite[Ch.~1]{m:yc} for a detailed description of the
algebraic structure of $\Y(\gl_N)$. The {\it reflection algebra\/}
$\Br(\gl_N)$ (see \cite{s:bc})
is isomorphic to the subalgebra of $\Y(\gl_N)$
generated by the coefficients of the series $s_{ij}(u)$
where the matrix $S(u)=[s_{ij}(u)]$ is given by
$S(u)=T(-u)^{-1}\ts T(u)$; see \cite{mr:rr}.

Now consider the universal enveloping algebra $\U(\gl_N)$
and let $E=[E_{ij}]$ be the $N\times N$ matrix whose
$(i,j)$ entry is the generator $E_{ij}$ of $\U(\gl_N)$.
For any complex parameter $\om\in\CC$
the mapping
\ben
T(u)\to \frac{u-E^{\tss t}-\om/4}{u-\om/4}
\een
defines an {\it evaluation homomorphism\/}
$\Y(\gl_N)\to\U(\gl_N)$; see e.g. \cite[Secs.~1.1-1.3]{m:yc}.
Its restriction to the reflection algebra is
a homomorphism $\Br(\gl_N)\to \U(\gl_N)$
given by
\beql{evra}
S(u)\mapsto \frac{u+\om/4}{u-\om/4}\cdot
\frac{u-E^{\tss t}-\om/4}{u+E^{\tss t}+\om/4}.
\eeq

Now, let $\la$ be a partition of $n$ and
suppose that $\Uc$ is a standard tableau
of shape $\la$.
Consider the corresponding idempotent $E_{\Uc}$
and the sequence of contents $(\si_1,\dots,\si_n)$
as defined in Corollary~\ref{cor:sym}.
The symmetric group $\Sym_n$ acts naturally
on the space \eqref{cntenpr}. By the Schur--Weyl duality
the subspace
\beql{sublugl}
L_{\Uc}=E_{\Uc}(\CC^N)^{\ot n}
\eeq
is nonzero if the number of boxes in the first
column of $\la$ does not exceed $N$;
in this case $L_{\Uc}$ is
an irreducible representation of $\gl_N$
with the highest weight $\la$.

\bpr\label{prop:invcogl}
The subspace $L_{\Uc}$
is invariant under the action of $\Br(\gl_N)$
given by
\begin{multline}
S(u)\mapsto \frac{u-\om/4}{u+\om/4}\ts
\ts R_{0n}(-u-\si_n-\om/4)^{-1}\ts
\dots R_{01}(-u-\si_1-\om/4)^{-1}\\
{}\times R_{01}(u-\si_1-\om/4)\ts
\dots R_{0n}(u-\si_n-\om/4).
\non
\end{multline}
Moreover, the representation of $\Br(\gl_N)$ on $L_{\Uc}$
is isomorphic to the evaluation module
$L_{\Uc}$ obtained via the homomorphism \eqref{evra}.
\epr

\bpf
Using the matrix notation and regarding the matrix $E$ as
the element
\ben
E=\sum_{i,j=1}^N e_{ij}\ot E_{ij}\in \End\CC^N\ot\U(\gl_N)
\een
we find that the image of the transposed matrix $E^{\tss t}$ in
\eqref{endspp} coincides with the element $X_0$ given by
\ben
X_0=P_{01}+P_{02}+\dots+P_{0n}.
\een
The argument is now competed in the same way as the proof
of Proposition~\ref{prop:invco} with the use of the image
of the identity \eqref{jmm} in the group algebra
of the symmetric group.
\epf


\begin{thebibliography}{99}

\bibitem{aacfr:rp}
{D. Arnaudon, J. Avan, N. Cramp\'e, L. Frappat and E. Ragoucy},
{\it $R$-matrix presentation for super-Yangians
$Y({\rm osp}(m\vert 2n))$},
{J. Math. Phys.}  {\bf 44}  (2003), 302--308.

\bibitem{amr:rp}
{D. Arnaudon, A. Molev and E. Ragoucy},
{\it On the $R$-matrix realization of Yangians
and their representations},
Annales Henri Poincar\'e {\bf 7} (2006), 1269--1325.

\bibitem{bw:bl}
{J. Birman and H. Wenzl},
{\it Braids, link polynomials and a new algebra},
Trans. AMS {\bf 313} (1989), 249--273.

\bibitem{b:aw}
R. Brauer,
{\it On algebras which are connected with
the semisimple continuous groups},
Ann. Math. {\bf 38} (1937), 854--872.

\bibitem{c:sb}
{I. V. Cherednik}, {\it On special bases of irreducible
finite-dimensional representations of the degenerate affine Hecke
algebra}, {Funct. Analysis Appl.} {\bf 20} (1986), 87--89.

\bibitem{c:ni}
{I. V. Cherednik},
{\it A new interpretation of Gelfand--Tzetlin bases}, {Duke Math. J.}
{\bf 54}
(1987),
563--577.

\bibitem{d:ha}
{V. G. Drinfeld},
{\it Hopf algebras and the
quantum Yang--Baxter equation},
{Soviet Math. Dokl.} {\bf
32} (1985), 254--258.

\bibitem{d:qg}
{V. G. Drinfeld},
{\it Quantum Groups},
in ``International Congress of Mathematicians (Berkeley, 1986)'',
Amer. Math. Soc.,
Providence RI,
1987,
pp. 798--820.

\bibitem{im:fp}
A. P. Isaev and A. I. Molev,
{\it Fusion procedure for the Brauer algebra},
Algebra i Analiz {\bf 22} (2010), 142--154.

\bibitem{imo:fp}
A. P. Isaev, A. I. Molev and O. V. Ogievetsky,
{\it Fusion procedure for the Birman--Murakami--Wenzl
algebra and evaluation homomorphisms for quantum
algebras}, in preparation.

\bibitem{imo:ih}
A. P. Isaev, A. I. Molev and A. F. Os'kin,
{\it On the idempotents of Hecke algebras},
Lett. Math. Phys. {\bf 85} (2008), 79--90.

\bibitem{io:bs}
A.~P.~Isaev and O.~V.~Ogievetsky,
{\it On Baxterized solutions of reflection
equation and integrable chain models},
Nuclear Phys.\  B {\bf 760} (2007), 167--183.

\bibitem{io:jm}
A.~P.~Isaev and O.~V.~Ogievetsky,
{\it Jucys-Murphy elements
for Birman--Murakami--Wenzl
algebras}, in: ``Proc. of Int. Workshop Supersymmetries
and Quantum Symmetries", Dubna 2009;
{\tt arXiv:0912.4010}.

\bibitem{j:yo}
{A. Jucys}, {\it On the Young operators of the symmetric group},
{Lietuvos Fizikos Rinkinys} {\bf 6} (1966), 163--180.

\bibitem{j:fy}
{A. Jucys}, {\it Factorization of Young projection operators for
the symmetric group}, {Lietuvos Fizikos Rinkinys}
{\bf 11} (1971), 5--10.

\bibitem{krs:yb}
{P. P. Kulish, N. Yu. Reshetikhin and E. K. Sklyanin},
{\it Yang--Baxter equation and representation theory},
{Lett. Math. Phys.}
{\bf 5}
(1981),
393--403.

\bibitem{lr:rh}
R. Leduc and A. Ram,
{\it A ribbon Hopf algebra approach to the
irreducible representations of
centralizer algebras: The Brauer, Birman-Wenzl and
type A Iwahori-Hecke algebras}, Adv. Math.
{\bf 125} (1997), 1--94.

\bibitem{m:fp}
{A. I. Molev},
{\it On the fusion procedure for the symmetric group},
Reports Math. Phys. {\bf 61} (2008), 181--188.

\bibitem{m:yc}
A. Molev,
{\it Yangians and classical Lie algebras}, Mathematical
Surveys and Monographs, 143. American Mathematical Society,
Providence, RI, 2007.

\bibitem{mr:rr}
{A. I. Molev and E. Ragoucy},
{\it Representations of reflection algebras},
Rev. Math. Phys. {\bf 14} (2002),
317--342.

\bibitem{m:nc}
{G. E. Murphy},
{\it A new construction
of Young's seminormal representation
of the symmetric groups }, {J. Algebra} {\bf 69} (1981), 287--297.

\bibitem{n:yc}
{M. Nazarov}, {\it Yangians and Capelli identities}, in:
``Kirillov's Seminar on Representation Theory" (G.~I.~Olshanski,
Ed.), {Amer. Math. Soc. Transl.} {\bf 181}, Amer. Math. Soc.,
Providence, RI, pp. 139--163 (1998).

\bibitem{n:yo}
{M. Nazarov},
{\it Young's orthogonal form for Brauer's
centralizer algebra},
J. Algebra {\bf 182} (1996), 664--693.

\bibitem{n:rt}
{M. Nazarov},
{\it Representations of twisted Yangians
associated with skew Young diagrams},
Selecta Math. (N.S.) {\bf 10} (2004), 71-129.

\bibitem{n:mh}
{M. Nazarov},
{\it A mixed hook-length formula for affine Hecke algebras},
European J. Combin. {\bf 25} (2004), 1345--1376.

\bibitem{o:ty}
{G. Olshanski},
{\it Twisted Yangians and infinite-dimensional classical Lie algebras},
in: ``Quantum Groups (Leningrad, 1990)'',
{Lecture Notes in Math.}
{\bf 1510},
Springer,
Berlin,
1992,
pp.~103--120.

\bibitem{rtf:ql}
{N. Yu. Reshetikhin, L. A. Takhtajan and L. D. Faddeev},
{\it Quantization of Lie Groups and Lie algebras}, {Leningrad Math. J.}
{\bf 1}
(1990),
193--225.

\bibitem{r:cs}
{H. Rui},
{\it A criterion on the semisimple Brauer algebras},
J. Comb. Theor., A {\bf 111} (2005), 78--88.

\bibitem{s:bc}
{E. K. Sklyanin},
{\it Boundary conditions for integrable quantum systems}, {J. Phys.}
{\bf A21}
(1988),
2375--2389.

\bibitem{w:sb}
H. Wenzl,
{\it On the structure of Brauer's centralizer algebras},
Ann. Math. (2) {\bf 128} (1988), 173--193.

\bibitem{w:cg}
{H. Weyl},
{\it Classical groups, their invariants and representations},
{Princeton Univ. Press},
Princeton NJ,
1946.

\bibitem{zz:rf}
{A. B. Zamolodchikov and Al. B. Zamolodchikov},
{\it Factorized $S$-matrices in two dimensions as the exact solutions
of certain relativistic quantum field models},
{Ann. Phys.} {\bf 120} (1979), 253--291.

\end{thebibliography}
\end{document}